\documentclass[11pt]{amsart}
\usepackage{amsmath,amssymb}
\usepackage{graphicx}
\usepackage{amsfonts,amscd,latexsym}

\renewcommand{\Hat}{\widehat}
\newcommand{\Ha}{{\Hat a}}
\newcommand{\Hb}{{\Hat b}}
\newcommand{\Hc}{{\Hat c}}
\newcommand{\He}{{\Hat e}}
\newcommand{\Hq}{{\Hat q}}

\newcommand{\oq}{{\otimes{\Q}}}

\newcommand{\Tf}{{\Tilde f}}

\newcommand{\Tev}{{\Tilde \ev}}
\newcommand{\Tka}{{\Tilde \ka}}

\newcommand{\TJ}{{\Tilde J}}
\newcommand{\TP}{{\Tilde P}}
\newcommand{\TQ}{{\Tilde Q}}
\newcommand{\im}{{\rm im}}
\newcommand{\QED}{{\hfill$\Box$\MS}}
\newcommand{\Vol}{{\rm Vol}}
\newcommand{\PD}{{\rm PD}}
\newcommand{\Aa}{{\mathcal A}}

\newcommand{\ev}{{\rm ev}}

\newcommand{\p}{{\partial}}
\newcommand{\Ee}{{\mathcal E}}

\newcommand{\Jj}{{\mathcal J}}

\newcommand{\Ppp}{{\mathcal P}}
\newcommand{\Mm}{{\mathcal M}}
\newcommand{\Cc}{{\mathcal C}}
\newcommand{\oMm}{{\overline {\Mm}}}
\newcommand{\ov}{\overline}
\newcommand{\be}{{\beta}}
\newcommand{\al}{{\alpha}}
\newcommand{\De}{{\Delta}}
\newcommand{\La}{{\Lambda}}

\newcommand{\io}{{\iota}}
\newcommand{\Si}{{\Sigma}}
\newcommand{\si}{{\sigma}}
\newcommand{\ga}{{\gamma}}

\newcommand{\la}{{\lambda}}
\newcommand{\eps}{{\epsilon}}
\newcommand{\Om}{{\Omega}}
\newcommand{\om}{{\omega}}
\newcommand{\de}{{\delta}}
\newcommand{\ka}{{\kappa}}
\newcommand{\Q}{{\mathbb Q}}
\newcommand{\R}{{\mathbb R}}
\newcommand{\C}{{\mathbb C}}

\newcommand{\T}{{\mathbb T}}
\newcommand{\SO}{{\rm SO}}

\newcommand{\Z}{{\mathbb Z}}
\newcommand{\CP}{{\mathbb CP}}
\newcommand{\ind}{{\rm ind }}
\newcommand{\Ham}{{\rm Ham}}
\newcommand{\Symp}{{\rm Symp}}

\newcommand{\Flux}{{\rm Flux}}
\newcommand{\Homeo}{{\rm Homeo}}

\newcommand{\Diff}{{\rm Diff}}

\newcommand{\Map}{{\rm Map}}
\newcommand{\vol}{{\rm vol}}
\newcommand{\Nn}{{\mathcal N}}

\renewcommand{\Tilde}{\widetilde}
\newcommand{\TM}{{\Tilde M}}
\newcommand{\TY}{{\Tilde Y}}
\newcommand{\TX}{{\Tilde X}}
\newcommand{\TSi}{{\Tilde \Si}}

\newcommand{\Bb}{{\mathcal B}}

\newcommand{\Tom}{{\Tilde\om}}
\newcommand{\TOm}{{\Tilde\Om}}

\newcommand{\Tp}{{\Tilde{p}}}
\newcommand{\Tq}{{\Tilde{q}}}
\newcommand{\Tg}{{\Tilde{g}}}
\newcommand{\Tx}{{\Tilde{x}}}

\newcommand{\labell}[1] {\label{#1}}

\newtheorem{theorem}{Theorem}[section]

\newtheorem{thm}[theorem]{Theorem}

\newtheorem{cor}[theorem]{Corollary}

\newtheorem{lemma}[theorem]{Lemma}

\newtheorem{prop}[theorem]{Proposition}

\newtheorem{rmk}[theorem]{Remark}
\newtheorem{example}[theorem]{Example}
\numberwithin{equation}{section}

\newcommand{\SSS}{{\smallskip}}
\newcommand{\MS}{{\medskip}}

\newcommand{\GW}{{\rm GW}}
\newcommand{\PGW}{{\rm PGW}}
\newcommand{\NI}{{\noindent}}

\begin{document}
\title{The symplectomorphism group of a blow up}
\author{Dusa McDuff}\thanks{partially supported by NSF grants DMS 0305939 and 0604769}
\address{Department of Mathematics,
 Stony Brook University, Stony Brook, 
NY 11794-3651, USA}
\email{dusa@math.sunysb.edu}
%\urladdr{http://www.math.sunysb.edu/\~{}dusa}
\date{October 4, 2006, with minor revisions May 29, 2007}

\begin{abstract}
We study the relation between the symplectomorphism group 
$\Symp\, M$ of a closed connected symplectic manifold $M$ and 
the symplectomorphism and diffeomorphism groups $\Symp\,\TM$ and $\Diff \,\TM$ of its one point blow up $\TM$. 
 There are three main arguments.  The first 
shows that for any oriented $M$ the
     natural map from $\pi_1(M)$ to $\pi_0(\Diff \TM)$ is often injective.    The second argument applies  
    when $M$ is simply connected and detects
      nontrivial elements in the homotopy group $\pi_1(\Diff \TM)$ 
    that persist into the space of
self homotopy equivalences of $\TM$. 
      Since it  uses purely homological arguments,
       it applies to
 $c$-symplectic manifolds $(M,a)$, that is, to manifolds of dimension $2n$ that support a class $a\in H^2(M;\R)$ such that $a^n\ne 0$.    
 The third argument uses the symplectic structure on $M$ and 
 detects nontrivial elements in the (higher) homology of 
 $B\Symp\, \TM$ using characteristic classes defined 
 by parametric Gromov--Witten invariants.
 Some results about many point blow ups are also obtained. 
 For example we show that if $M$ is the $4$-torus with
 $k$-fold blow up $\TM_k$ (where $k>0$) then
   $\pi_1(\Diff \TM_k)$ is not generated by the groups
   $\pi_1\big(\Symp\, (\TM_k, \Tom)\bigr)$ as 
   $\Tom$ ranges over the set of all symplectic forms on $\TM_k$. 
\end{abstract}

\maketitle
%\begin{center} Preliminary version\end{center}

  %%%%%%%%%%%%%%%%%%%%%%%%%%%%%%%%%%%%%%%%%%%%%%%%%%%%%%%%%%%%%%%%%
  \section{Main results}
  %%%%%%%%%%%%%%%%%%%%%%%%%%%%%%%%%%%%%%%%%%%%%%%%%%%%%%%%%%%%%%%%%

In this note we investigate the relation between the symplectomorphism group of a closed connected symplectic manifold $(M,\om)$ of dimension $2n\ge 4$ and of its one point blow up $(\TM,\Tom_\eps)$.  (Here $\eps>0$ measures the size of the blow up.)  When $\dim M = 4$ and $\eps$ is sufficiently small, it is known from work of Lalonde--Pinsonnault~\cite{LP} 
and K\c edra~\cite{K2} that in several cases there is a homotopy equivalence between $\Symp\, \TM: = \Symp(\TM,\Tom_\eps)$ and
 $\Symp(M,p): = \Symp(M,p,\om)$, 
  the group of symplectomorphisms of $M$ that fix the point $p$. Our aim here is to understand 
 both the role played by the symplectic form and what happens in higher dimensions.  Throughout we investigate stable phenomena that persist in the diffeomorphism group $\Diff$ 
 and under various kinds of deformation; see Remark~\ref{rmk:XX}.
  
     We will present three different arguments.  The first studies
     a natural map $\Tf_*:\pi_1(M)\to \pi_0(\Diff \TM)$ defined for any  oriented manifold $M$, and shows that $\Tf_*$
      is often injective.    
      The second works best for simply connected   $M$ and
     detects
      nontrivial elements in the homotopy group $\pi_1(\Diff \TM)$ that persist into $\pi_1(\TM^{\TM})$,
       where $X^X$ denotes the space of continuous selfmaps of $X$.
      Since it  uses purely homological ideas,
       this argument applies to
 $c$-symplectic manifolds $(M,a)$, that is, to manifolds of dimension $2n$ that support a class $a\in H^2(M;\R)$ such that $a^n\ne 0$.  
 In \S\ref{ss:4} we add some symplectic ideas to investigate 
 the difference between $\pi_1(\Diff \TM_k)$
 and $\pi_1(\Symp\, \TM_k)$ for a $k$-fold blow up of a $4$-manifold. The third argument, presented in \S\ref{ss:homol},
  needs the symplectic structure on $M$ and 
 detects nontrivial elements in the (higher) homology of $B\Symp\,
 \TM$.    It uses 
a characteristic class on $B\Symp\,\TM$
defined by a suitable parametric Gromov--Witten invariant.

Before explaining our results in more detail, we need some preparation.
Recall that the classifying space $B\Symp(M,p)$ is the total space of the universal $M$-bundle over $B\Symp\, M$
\begin{equation}\labell{eq:unib}
  M\to B\Symp(M,p) \stackrel{\pi}\to B\Symp\, M.
\end{equation}
  (Later we sometimes denote this bundle by $M_{\Symp}\to B\Symp\,M$.)
The pullback of this bundle over $\pi: 
  B\Symp(M,p) \to B\Symp\, M$ is the universal $M$-bundle with section.
  In other words, a map $\phi:X\to B\Symp(M,p)$ determines and is determined up to homotopy by a triple $(P_\phi, \pi_\phi, s)$, where
$\pi_\phi: P_\phi\to X$ is a symplectic $M$-bundle and $s:X\to P_\phi$ is a section. Blowing up along this section gives a map\footnote
{
We take the target space to be $B\Diff$ rather than $B\Symp$ since then we can use the pointwise blow up of algebraic geometry. Because $B\Symp(M,p)$ is 
noncompact it is not clear that there is a symplectic blow up of uniform size $\eps>0$
over the whole of $B\Symp(M,p)$.}
$$
\be:B\Symp(M,p)\to B\Diff\TM.
$$
The composite map 
\begin{equation}\labell{eq:f}
\Tf: M\to B\Diff \TM
\end{equation}
represents the $\TM$-bundle over $M$ that is obtained from the trivial bundle $M\times M\to M$ by blowing up along the diagonal.

More generally, let $M$ be an oriented even dimensional smooth manifold,
choose a Hermitian metric
%%DD unitary framing of the tangent space $T_pM$ 
on  the tangent space $T_pM$ at $p$ that is compatible with the orientation (and with the symplectic form if one is given), and
denote by $\Diff^{U}(M,p)$ the subgroup of $\Diff(M,p)$ consisting of diffeomorphisms whose derivative at $p$ preserves the Hermitian metric and so is unitary.
%%this framing.
Then, maps $X\to B\Diff^{U}(M,p)$ classify smooth $M$-bundles over $X$ provided with sections $\si$ that have Hermitian normal bundles. 
Hence again one can blow up along the section to get a map
  $$
  \be^U:B\Diff^{U}(M,p)\to B\Diff\TM.
  $$
  Because a symplectic vector bundle has a complex structure
   that is unique up to homotopy, 
for symplectic $M$ there is a commutative diagram
$$
\begin{array}{ccc} B\Symp^{U}(M,p)&\stackrel{\simeq}\to& B\Symp(M,p)\\
\io\downarrow && \be\downarrow\\
B\Diff^{U}(M,p)&\stackrel{\be^U}\to& B\Diff\TM
\end{array}
$$
in which the inclusion on the top line is a homotopy equivalence.

Suppose now that $M$ is an oriented manifold of dimension $2n$,
and consider the composite
$$
\Tf_*:\pi_1M\to \pi_1\bigl(B\Diff(M,p)\bigr)\stackrel{\be^U}\to \pi_1(B\Diff\TM).
$$
 (The first map here is induced by (\ref{eq:unib}), while the second  is well defined  because 
an oriented $\R^{2n}$ bundle over $S^1$ has a unique 
homotopy class of complex structures.)  Because there is a fibration
$$
\Diff\,M\stackrel{\ev}\to M\to B\Diff(M,p),\quad \ev(\phi) :=\phi(p),
$$
the kernel of ${\Tf_*}$ contains the evaluation subgroup $\ev_*\bigl(\pi_1(\Diff\,M)\bigr)$.

\begin{thm}\labell{thm:pi0} For any closed oriented $2n$-manifold 
$M$  the kernel of 
$$
{\Tf_*}: \pi_1M\to \pi_1(B\Diff\,\TM)
$$
 acts trivially on 
$\pi_iM$ for $1\le i<\dim M$.  In particular, it is 
contained in the center of $\pi_1M$ and vanishes if $M$ itself is a blow up. It also vanishes if the Euler characteristic $\chi(M)$ is nonzero.
\end{thm}

\begin{rmk}\labell{rmk:ev}\rm  We will see in \S2 that the kernel 
 of $\Tf_*$ is contained in the image of the evaluation map
 $\Tev_*: \pi_1(M^{\TM})\to \pi_1M$ given by $\Tev(\phi)=
 \phi(\Tp)$, where $\Tp\in \TM$ is fixed. (Here $M^{\TM}$ denotes the space of maps $\TM\to M$, with base point at the blow down $\phi_0:\TM\to M$.)
   This evaluation map is very similar to the usual one 
 $\ev_*:\pi_1(M^M)\to \pi_1M$ considered by Gottlieb~\cite{Gott}, though it might in principle have a larger image.  The properties mentioned above are some of those that
 Gottlieb established for $\im\,\ev_*$. 
 \end{rmk}
% may not equal the whole of 
%${\rm im }\ev_H$, but
%it does contain $\ev_H(\pi_1(\Diff M))$ because of the exactness of the sequence 
%$$
%\pi_1(\Diff M)\stackrel{\ev_H} \to \pi_1(M)\to \pi_1(B \Diff M).
%$$
%Note also that ${\rm im }\ev_H$ often vanishes.
%  According to Gottlieb~\cite{}, it is contained in the center of $\pi_1(M)$ and also acts trivially on $\pi_*(M)$.  It coincides with the center
%if $M$ is a $K(\pi,1)$, and vanishes when the Euler characteristic of $M$ is nonzero.

We now consider $\pi_2$.
Denote by 
\begin{equation}\labell{eq:Ic}
I_c:\quad \pi_2(B\Diff^{U}(M,p))\to \pi_2(BU(n)) \cong \Z. 
 \end{equation}
the homomorphism induced by 
taking the derivative at $p$. In the case when $[\si]\in
\pi_2(B\Diff^{U}(M,p))$ is in the image of $B\Symp(M,p)$ and so corresponds to a section $s_\si$
 of a symplectic bundle $P\to S^2$ with its canonical Hermitian structure, then $I_c([\si])$ is the value on $s_\si$ of the  first Chern class $c_1^V$ of the vertical tangent bundle.

 If $(M,a)$ is $c$-symplectic and simply connected, there is another useful homomorphism
$$
 I_a: \pi_2(B\Diff^{U}(M,p))\to \R
$$
 that we now define.
Consider the universal $M$-bundle 
\begin{equation}\labell{eq:Mb}
M\to B\bigl(\Diff(M,p)\cap \Diff_0M\bigr)\stackrel{\pi} \to B\Diff_0M,
\end{equation}
where $\Diff_0$ denotes  the identity component  of $\Diff$.
  Denote by $\Cc: = \Cc_M$ the (open) cone 
  consisting of all classes $a\in H^2(M;\R)$ such that $a^n\ne 0$.

\begin{lemma}[\cite{KMc}]\labell{le:coup}
If $H^1(M;\R)=0$ the restriction map 
$$
H^2\Bigl(B\bigl(\Diff(M,p)\cap \Diff_0M\bigr);\R\Bigr)\to H^2(M;\R)
$$
 is surjective.  
Moreover, each $a\in \Cc$ has a unique extension to a class
%%D Ta to Ha
$\Ha$ such that 
the
fiberwise integral $\pi_!(\Ha): = \int_M {\Ha}^{n+1}\in H^2(B\Diff_0M)$ vanishes.
\end{lemma}
\NI {\it Sketch of proof.}\,\, Consider the Leray--Serre spectral sequence for 
the fibration (\ref{eq:Mb}). Since $E_2^{2,1} = 0$ and $\Cc$ is open, the first statement will hold provided that 
the differential $d:=d_{0,2}^3$ vanishes on all $a\in \Cc$.   But $0 = d(a^{n+1}) = (n+1) a^n da$ which implies that $da=0$ since multiplication by $a^n$ induces an isomorphism $E_2^{3,0}\to E_2^{3,2n}$.  If $\Ha'$ is any extension of $a$, it is easy to check that
$$
\Ha: = \Ha' - \frac 1{\la (n+1)} \pi^* \bigl(\pi_!(\Ha')^{n+1}\bigr), 
$$
satisfies the conditions if $\la:= \int_M a^n$. (Here $\pi_!:H^{2n+2}\to H^2$ 
denotes the cohomology push forward, i.e. integration over the fiber.)
\QED

The pullback of $\Ha$ to the total space of an $M$-bundle $P\to Z$ is often
 called  the {\bf coupling class}; we will also call it the {\bf normalized extension} of $a$. On a product bundle,
for example, it  is just the obvious pullback $pr_M^*(a)$ of $a$.
Since each $[\si] \in \pi_2(B\Diff^{U}(M,p))$ has a unique lift to an element (also called $[\si]$) of 
 $\pi_2(B(\Diff^{U}(M,p)\cap \Diff_0))$, we may define $I_a$ by setting:
\begin{equation}\labell{eq:Ia}
 I_a([\si]): = \int_{S^2} \si^*(\Ha).
\end{equation}
If we think of $[\si]$ as an $M$-bundle $P\to S^2$  with section $s_\si$, then $ I_a([\si])$ is given by evaluating 
the coupling class $\Ha$ on the section $s_\si$.

  We say that $[\si]\in \pi_2(B\Diff^{U}(M,p))$ is 
  {\bf homologically visible} if it has nonzero image under the 
  homomorphism
  $I_c\oplus I_a$ for some $a\in \Cc$.
  
  \begin{thm}\labell{thm:csymp}  Let $M$ be a simply connected 
  $c$-symplectic manifold.  Then every homologically visible
  class $[\si]\in \pi_2(B\Diff^{U}(M,p))$ has nonzero image under
  the composite map  
  $ \pi_2(B\Diff^{U}(M,p)) \to \pi_2\bigl(B\Diff \TM\bigr)\to 
  \pi_2\bigl(B (\TM^{\TM})\bigr).$
     \end{thm}
    
\begin{cor}\labell{cor:csymp} Let $(M,a)$ be a $c$-symplectic manifold.  Then there is a  homomorphism $\Tf_*: \Z \oplus \pi_2(M) \to \pi_2\bigl(B\Diff\TM\bigr)$ whose kernel is 
contained in the kernel of the rational Hurewicz homomorphism
$\pi_2(M)\to H_2(M;\Q)$.
\end{cor}
    
\begin{proof}
First assume that   $\pi_1(M) = 0$. 
The boundary map $\pi_{*+1}(B\Diff_0)\to 
    \pi_{*}(M)$ in the fibration~(\ref{eq:Mb}) desuspends to
    the map $\ev_*: \pi_*(\Diff_0)\to  \pi_*(M)$ induced by 
    evaluation 
    $$
    \ev: \Diff_0 \to M,\quad \phi\mapsto \phi(p).
    $$
      Since $\ev$ induces the trivial map on $\pi_2\oq$, 
      the kernel of the induced map 
    $f_*: \pi_2(M)\to \pi_2(B\Diff(M,p))$ consists of torsion elements.  This map $f_*$ can be described explicitly as follows.
    Given $\si:S^2\to M$ the element $f_*[\si]$ corresponds to the trivial bundle $M\times S^2\to S^2$ with section $gr_\si: z\mapsto (\si(z),z)$.  Since the normalized extension of the class $a\in \Cc_M$ to $M\times S^2$ is just $pr_M^*(a)$, we find that $I_a(f_*[\si]) = 
    a([\si])$.  Therefore,  because $\Cc$ is open,
the elements in $\im f_*$ that are not
 detected by some homomorphism $I_a$ have the form $f_*[\si]$ where 
$[\si]$ maps to zero  in $H_2(M;\Q)$.
    Since $\pi_2(SO(2n)/U(n)) \cong\Z$, each element in
     $\pi_2(B\Diff(M,p))$ has $\Z$ different lifts to
     $\pi_2(B\Diff^{U}(M,p))$, distinguished by the values of $I_c$.
    Now apply Theorem~\ref{thm:csymp}.
     The general case is proved in \S3.  \end{proof}

\begin{rmk}\rm 
When $\pi_1(M)=0$, the above argument shows that the elements persist in $\pi_2\bigl(B (\TM^{\TM})\bigr)$.
However this may not hold in general.
\end{rmk}

When $(M, \om)$ is symplectic, the map $\Tf: M\to  
B\Diff\TM$ defined in (\ref{eq:f}) factors through 
$B\Symp(\TM,\Tom_\eps)$ for suitably small $\eps$.
One can construct this lift
   explicitly as follows. 
Consider the product $M\times M$ with symplectic form $\Om: = \om
\times \om$. Then the diagonal is a symplectic submanifold, and so one can blow up normal to it by some amount $\eps_0>0$ to get a symplectic fibration
\begin{equation}\labell{eq:MM}
(\TM, \Tom_{\eps_0}) \to (\TP_{\eps_0}, \TOm_{\eps_0}) \to M.
\end{equation}
By shrinking the size of 
the blow up, one gets a corresponding fibration $(\TP_{\eps},
\TOm_{\eps})  \to M$ for 
every $\eps\le \eps_0$.  Note that the initial fibration
$(\TP_{\eps_0}, \TOm_{\eps_0}) \to M$ might depend on the choice 
of blow up.  
 However any two choices 
become equivalent
for small enough $\eps>0$.
Define 
\begin{equation}\labell{eq:fe}
\Tf_\eps: M\to B\Symp(\TM,\Tom_\eps)
\end{equation}
to be the classifying map of this fibration. 
Then for each $\si:S^2\to M$ the class $(\Tf_\eps)_*([\si])$ is represented by the pullback of $\TP_\eps\to M$ over $\si$.   
(K\c edra gives another description of this map in \cite{K2}.)
%If $\pi_1(M)=0$, this map lifts  to $B\Symp_0 = B\Ham$ vv
%%%
%\begin{equation}\labell{eq:fe0}
%\Tf_\eps: M\to B\Symp_0(\TM,\Tom_\eps)
%\end{equation}
%%%

Theorem~\ref{thm:csymp} implies that the subgroup $\Tf_*(\Z\oplus 0)$ is disjoint from the image of $(\Tf_\eps)_*(\pi_2(M))$ in $\pi_2(B\Diff\TM)$. It seems very likely that for any form $\Tom_\eps$ there are no elements in   $\pi_2\bigl(B\Symp(\TM, \Tom_\eps)\bigr)$ that map into $\Tf_*(\Z\oplus 0)$. 
This is true when $M=\C P^2$ since in this case  $\pi_2\bigl(B\Symp(\TM, \Tom_\eps)\bigr)$ is known to be $\Z$.
It is not clear how to prove this in general, though \S6 makes some progress  in the $4$-dimensional case: see 
Proposition~\ref{prop:diff2} and Corollary~\ref{cor:T4}.

The next result is based on detailed knowledge of the symplectic geometry of  the one point blow up $\C  P^2\#{\ov{\C P}}\,\!^2=:X$ of $\C P^2$.  In particular, it uses
the fact
(proved in  \cite{LM2}) that every symplectic form on 
  $X$ is diffeomorphic to some $\Tom_\eps$.

\begin{prop}\labell{prop:diff}  Let $X= \C  P^2\#
  {\ov{\C P}}\,\!^2$. 
  Then $\pi_2(B\Diff X)$ is not generated by the images of 
  $\pi_2(B\Symp(X,\om))$, as $\om$ varies over the space of all symplectic forms on $X$.
\end{prop}

\begin{rmk}\labell{rmk:X}\rm (i) During the 
 proof of this proposition we show
 that
 the obvious inclusion $S^2:=\C P^1\hookrightarrow \C P^2$
 gives rise not to the generator of 
$\pi_2(B\Symp X)$ (as one might initially expect) but to three times this generator.
\SSS

\NI (ii)  There are many blow up manifolds  that do not admit circle actions.   For example, in dimension $4$ Baldridge~\cite{Bald}
 has shown that if $X$ has $b^+>1$ and admits a circle action with a fixed point then its Seiberg--Witten invariants must vanish.  Since 
 manifolds that admit fixed point free circle actions must have zero Euler characteristic, this implies that no simply connected K\"ahler surface with $b^+>1$  admits a circle action.  On the other hand, 
Corollary~\ref{cor:csymp} implies that if $X$ is a blow up
and $\pi_1X=0$
 the groups $\pi_1(\Symp)$ and 
$\pi_1(\Diff)$ do not vanish.  Thus, as
 pointed out in K\c edra~\cite{K1}, there are plenty of examples of
 nonzero elements of $\pi_1(\Symp)$ or 
$\pi_1(\Diff)$  that cannot be represented by
a circle action.\footnote
{
By using the geometric decomposition of a Hamiltonian $S^1$ manifold provided by the moment map, it is easy to show  when $\dim M = 4$ that certain elements 
$\eta$ in $\pi_1(\Ham)$ cannot be represented by an $S^1$ action.  For example, if $(M,\om) = (S^2\times S^2, \la pr_1^*\si + pr_2^*\si)$, where $\la> 1$, one can take $\eta: = \eta_1+\eta_2$, where $\eta_1\in \pi_1\bigl(SO(3)\times \{id\}\bigr)\subset \pi_1\bigl(SO(3)\times SO(3)\bigr)$ rotates the first sphere and $\eta_2$ is one of the circle actions of infinite order. Since $\eta_1$ has finite order, the rational 
 Samelson product $[\eta,\eta]_{\Q}$ vanishes.   Buse in~\cite[Prop~2.1]{Bu2} gives an example of an element $\eta$ in $\pi_1(\Ham(\T^2\times S^2))$ for which $[\eta,\eta]_{\Q}$ does not vanish.
}  
One of the motivations of the present work was to understand the limits of K\c edra's argument.\end{rmk}

Our final  argument detects parts of the image of $H_*(M)$ in the homology of  $B\Symp\,\TM$.
%For simplicity, we will begin by assuming  that
If $\pi_1M = 0$,  
 the classifying map $\Tf_\eps: M\to B\Symp\,\TM$ of the fibration 
(\ref{eq:MM}) factors through the simply connected space
 $B\Symp_0\TM = B\Ham\,\TM$, where 
$\Ham$ denotes the Hamiltonian subgroup
$\Ham\subseteq \Symp_0$.  
We denote by $\Cc^*$ the subring
 of $H^*(M;\Q)$ generated\footnote
{Since the cone $\Cc\subset H^2(M)$ is open, this is the same as the subring generated by $\Cc$.}
by the classes in $H^2(M)$, 
and  define
$$
\Cc_*: = \{\ga\in H_*(M;\Q): \ka(\ga) = 0 \mbox{ for all } \ka\in \Cc^*\}.
$$
 
\begin{thm}\labell{thm:homol} Suppose that $\pi_1(M)=0$.  Then,  
the kernel of
$$
(\Tf_\eps)_*: H_*(M;\Q)\to H_*\bigl(B\Ham(\TM,\Tom_\eps);\Q\bigr)
$$
is contained in $\Cc_*$.
\end{thm}

When  $\pi_1(M)\ne 0$,  
Theorem~\ref{thm:pi0} implies that 
the bundle (\ref{eq:MM}) is usually not trivial over its $1$-skeleton. Thus its classifying map need not lift to
to $B\Ham$ (or even to $B\Symp_0$).  This is important because we detect the image of $(\Tf_\eps)_*$ by using characteristic classes, and
typically these live on $B\Ham$ rather than $B\Symp_0$ or $B\Symp$; cf. \S\ref{ss:homol}.
One way one might try to deal with this is to note that
because the fiberwise symplectic forms on (\ref{eq:MM}) 
 extend to the closed form $\TOm_\eps$, one can lift  the classifying map to the space $B\Ham^s$ defined in \cite{Mce}. One can then try to use characteristic classes coming from the coupling class on $M_{\Ham^s}$.  The problem here is
that when $H^1(M;\R)\ne 0$  the normalization condition in Lemma~\ref{le:coup} does not uniquely determine the extension $\Ha$ of $a=[\om]$ to the total space of an $M$-bundle $P\to Z$ unless 
$H^1(Z;\R) = 0$. 

We shall adopt an easier approach, restricting consideration to classes $[\si]\in H_d(M;\Q)$ that are represented by smooth maps $\si:Z\to M$ where $Z$ is a closed and {\it simply connected}  oriented 
$d$-manifold.  (For short we will call these smooth simply connected 
cycles.  If $\dim M> 4$ these generate
 the image of the homology of the universal cover of $M$; see
Lemma~\ref{le:sc}.)  The corresponding $\TM$-bundles $\TP\to Z$ 
are then classified by  maps $Z\to B\Symp_0\,\TM$.  These lift further to $B\Ham\,\TM$ because of the existence of the closed extension $\TOm_\eps$ of the fiberwise symplectic forms.  (This is another place where it is important that $\pi_1(Z)=0$; see \cite[Erratum]{LM}.)

Even with these conditions
the question of which classes  can be detected in $B\Ham\,\TM$ is quite subtle; see Example~\ref{ex:tor}.
Recall that
$(M,\om)$ is said to have the {\bf hard Lefschetz property}
if for all $k\le n$, the map $\cup [\om]^{k}: H^{n-k}(M;\R)\to 
 H^{n+k}(M;\R)$ is an isomorphism.  

 \begin{prop}\labell{prop:homol2} 
 Let $\si:Z\to M$ be a smooth and simply connected cycle representing a nonzero
 homology class  $[\si]\in H_{d}(M;\Q)$. 
 Let  $\Tf_\eps:Z\to B\Ham (\TM,\Tom_\eps)$ classify  
 the bundle
 $\TP_\eps\to Z$ obtained by blowing up $M\times Z\to Z$ along $gr_\si$.
 Then $(\Tf_\eps)_*[\si]\ne0$ if one of the following conditions holds.
 \smallskip
 
 \NI{\rm (i)} There is $\ka\in \Cc^*$ such that $\ka([\si])\ne 0$.
   
 \smallskip
 
 \NI{\rm (ii)}  
$Z=S^{2k}$ and $(M,\om)$ has the hard Lefschetz property.
\end{prop}

\begin{rmk}\labell{rmk:XX}\rm (i)  When we work in the symplectic category,
 we cannot control the size of the blow ups
and so must assume that they are arbitrarily small. We can say something about size only
when we make several blow ups; for example,
 if we blow up twice, the second blow up may have to be smaller than the first.  (For more precise statements see 
Propositions~\ref{prop:size} and~\ref{prop:kblow} and Corollary~\ref{cor:two}.)
\SSS

\NI
(ii) Our methods
 detect very stable symplectic phenomena that persist even in the diffeomorphism group.  
This should be contrasted with the work of 
 Seidel~\cite{Sei} who in the $4$-dimensional case
detects nonzero elements in $\pi_0(\Symp\,(X,\om_X))$ that disappear in $\pi_0(\Diff X)$. In his Example~1.13,  $(X,\om_X)$ is the monotone 
blow up of $\C P^2$ at $5$ points in which
 each point is blown up  to exactly one third of the size of the line
  $L$ in $\C P^2$. If one blows up each point slightly less (but still 
  keeps all the blow ups the same size) then at least some of these elements in $\pi_0(\Symp(X,\om_X))$ become zero:  for example, the class
  $L-E_1-E_2-E_3$ is represented by a Lagrangian sphere in the 
  monotone case, but by a symplectic sphere when the size is reduced.
 Therefore  the square of the Dehn twist in this sphere becomes trivial.  If the sizes of the blow ups are all different, then
 the above argument applies to a set of generators for 
 $\pi_0(\Symp^H(X,\om_X))$
so that this whole group vanishes. (Here
$\Symp^H$ denotes the subgroup acting trivially on $H_*(X)$.) 

  It would be 
 interesting to investigate the map  $$\pi_0(\Symp(X,\om_X))\to
 \pi_0(\Symp(X\times S^2,\om_X+\om_{S^2}))$$
 induced by $\phi\mapsto \phi\times id$ to see if it has nontrivial image.
  In contrast, the elements of $\pi_*(B\Diff\TM)$ considered in this paper do 
 remain nontrivial when $\TM$ is
 stabilized to 
a product $\TM\times Y$.\SSS

\NI (iii)  There is still rather little work that compares the
homotopy of  $\Symp$ and $\Diff$, or that discusses the dependence of $\pi_k(\Symp)$ on the symplectic form for manifolds of dimension $>4$.
Here are two notable exceptions.  Ruan~\cite{Ru}  
gives examples of $6$-manifolds such that
$\pi_0(\Symp_d)$ does not surject onto $\pi_0(\Diff)$, where
$\Symp_d$ is the group of diffeomorphisms that preserve that symplectic form up to deformation, while Seidel~\cite{Sei2}
studies the changes in 
$\pi_k\bigl(\Symp(\C P^n\times \C P^m,\om)\bigr)$ as $[\om]$ varies.
\end{rmk}

\NI
{\bf Acknowledgements.}  I thank Jarek K\c edra for useful comments
on an earlier draft of this paper, Daniel Gottlieb for help in 
understanding the evaluation subgroup and the referee 
for pointing out various small inconsistencies.

%%%%%%%%%%%%%%%%%%%%%%%%%%%%%%%%%%%%%%%%%%%%%%%%%%%%%%%%%%%%%
 \section{The fundamental group}
%%%%%%%%%%%%%%%%%%%%%%%%%%%%%%%%%%%%%%%%%%%%%%%%%%%%%%%%%%%%%

In this section we first discuss some homotopy theoretic approaches to the 
above questions, and then prove Theorem~\ref{thm:pi0}.
Let $M$ be a closed oriented manifold and consider the diagram
$$
\begin{array}{ccccccc}
\pi_k(\Diff \TM)&\longrightarrow&\pi_k(\Ee/U)&\longrightarrow &\pi_k\bigl(B\Diff^U(\TM,\Si)\bigr)&\longrightarrow&\pi_k(B\Diff\TM)\\
&&&&\simeq\uparrow &&\\
\pi_k(\Diff M)&\stackrel{\ev}\longrightarrow&
\pi_k(J\!M)&\longrightarrow&\pi_k\bigl(B\Diff^U(M,p)\bigr)
&\longrightarrow&\pi_k(B\Diff M).
\end{array}
$$
Here $J\!M$ denotes the space of pairs $(x,J_x)$, where $x\in M$ and $J_x:T_xM\to T_xM$ is an almost complex structure compatible with the orientation of $M$.  Thus there is a Hurewicz fibration 
$SO(2n)/U(n)\to J\!M\to M$.  Since $\Diff M$ acts transitively 
on $J\!M$  the bottom row in the above diagram
is exact. The top row is also exact. Here we identify the exceptional divisor $\Si$ with
$\C P^{n-1}$ and denote by $\Ee$ the space 
of all smooth embeddings of $\Si$ into $                                                                                                                   \TM$ that extend to diffeomorphisms of $\TM$. Thus an element of $\Ee$ factors (nonuniquely)  as 
$$
\C P^{n-1}\equiv \Si\stackrel f\to \Si\hookrightarrow \TM\stackrel{\Tg}\to \TM.
$$
If $\Ee/U$ is the quotient of $\Ee$
by the action of $U(n)$ on the domain, $\Diff(\TM)$ acts transitively on $\Ee/U$ with kernel 
 $\Diff^U(\TM,\Si)$ equal to
  the stabilizer of the inclusion map $\Si\hookrightarrow \TM$.   By standard arguments the groups $\Diff^U(M,p)$ and
$\Diff^U(\TM,\Si)$ are homotopy equivalent to the subgroups where the diffeomorphisms are unitary in a neighborhood of $p$ and $\Si$.  Thus these two groups are homotopy equivalent. 

We are interested in understanding the zig-zag composite
$$
\pi_k(J\!M)\longrightarrow
\pi_k\bigl(B\Diff^U(M,p)\bigr)\stackrel{\simeq}\longrightarrow \pi_k\bigl(B\Diff^U(\TM,\Si)\bigr)\longrightarrow
 \pi_k\bigl(B\Diff \TM\bigr).
$$
The first map is relatively understandable, 
since the previous discussion shows that its kernel
is closely related to the image of the evaluation map
$$
\ev: \Diff M\to M,\quad g\mapsto g(p).
$$
Therefore, the problem is to understand the kernel of the last map, i.e. the image of $\pi_k(\Ee/U)\to 
\pi_k\bigl(B\Diff^U(\TM,\Si)\bigr)$.  In the  $4$-dimensional 
symplectic situation considered by K\c edra, the space analogous to 
$\Ee/U$ consists of symplectically embedded $2$-spheres and is often 
contractible.\footnote
{
Suppose for example,  that $[\om]$ is integral and $\eps $ has the form $1/N$ for some integer $N$.  Then 
  the exceptional divisor $\Si$ is a curve of minimal 
  energy in $(\TM, \Tom_\eps)$ and so has a unique embedded $J$-holomorphic representative for each $\om$-tame $J$.  This implies that there is a homotopy equivalence from the contractible space $\Jj(\om)$ of $\om$-tame $J$ to $\Ee/U$: see Lalonde--Pinsonnault \cite[\S2]{LP}.
}
However, this need not hold in the smooth case
or in higher dimensions.

We now restrict to the case $k=1$. We will prove 
Theorem~\ref{thm:pi0} 
by showing that the kernel of $\pi_1 M\to \pi_1(B\Diff \TM)$ is the image of $\pi_1(\Bb)$, where $\Bb\subset M^{\TM}$ is  the space of smooth blow down maps.

To this end, fix a blow down map $\phi_0: (\TM,\Si) \to (M,p)$ and denote by $\Bb$ the set of 
%\footnote
%{
%One could ask that $\Tg$ be smooth, but must allow $h$ to be nonsmooth at the base point $p$.}
smooth maps conjugate to $\phi_0$. More precisely 
 $$
 \Bb: = \{h\circ \phi_0\circ  \Tg\in \Map^\infty(\TM,M): h\in \Homeo_0 M, \Tg\in \Diff_0\TM\},
 $$
 where $\Homeo M$ is the group of homeomorphisms\footnote
 { 
It is convenient to relax the smoothness conditions on $h$ since then
 for any diffeomorphism $\Tg: (\TM,\Si)\to (\TM,\Si)$ there is
 $h$ satisfying $h\circ\phi_0 = \phi_0\circ\Tg$. However, because
 $h\circ \phi_0\circ  \Tg$ and $\Tg$ are both smooth, $h$ is smooth except at the basepoint $p$.
 }
  of $M$ and the subscript $0$ denotes that we consider the identity component.  Thus $\phi:\TM\to M$ is an element of $\Bb$ iff 
  \begin{itemize} \item there
  is $q\in M$  such that $\phi: \TM\setminus\{\phi^{-1}(q)\} \to M\setminus q$ is a diffeomorphism, and\SSS
  
  \item  $\phi^{-1}(q) \in \ov\Ee_0$, where $\ov\Ee_0$ is the 
  identity component of 
   $\ov\Ee: = %\Ee/(\Diff \Si)= 
   (\Diff \TM)|_{\Si}/(\Diff\Si)$, i.e. the space 
 of smoothly embedded 
submanifolds isotopic to $\Si$. 
    \end{itemize}
(For then there is $\Tg\in \Diff_0\TM$ such that $\Tg(\phi^{-1}(q)) = \Si$ so that $\phi$ can be written as $h\circ\phi_0\circ\Tg$ for suitable $h$.)  
  
Although the elements $\phi\in \Bb$ do not have unique factorizations 
into pairs $(h,\Tg)$, each $\phi$
 induces a diffeomorphism $\TM\setminus\Tg^{-1}(\Si)\to M\setminus \{h(p)\}$.
In particular, there are well defined maps 
 $$
 \ev_B: \Bb\to M,\quad h\circ \phi_0\circ  \Tg \mapsto h(p),
 $$
and
 $$
 F: \Bb\to \ov\Ee_0,\quad h\circ \phi_0\circ  \Tg \mapsto \Tg^{-1}(\Si).
 $$
Observe also that $h\circ \phi_0\circ  \Tg= \phi_0$ 
iff $h = \Tg^{-1}$ on 
 $M\setminus \{p\} \equiv \TM\setminus \Si$.
 
Theorem~\ref{thm:pi0} follows immediately from the next two results.
As in Remark~\ref{rmk:ev},
$\Tev_*: \pi_1(M^{\TM})\to \pi_1M$ denotes the evaluation map given by $\Tev(\phi)=
 \phi(\Tp)$, where $\Tp\in \TM$ is fixed.   Note that this does 
 {\it not} restrict to $\ev_B$ on $\Bb$ since the base point $\Tp$ need not lie in
 $\Tg^{-1}(\Si_0)$.  Nevertheless, the next result shows that the two maps have comparable effects on $\pi_1(\Bb)$.

%denote by $\Tev:M^{\TM}\to M$ the evaluation map at $\Tp\in \Si$ 
%from the space $ M^{\TM}$ of continuous maps $\TM\to M$, where 
%the  basepoint of $ M^{\TM}$ is the blow down map $\phi_0:\TM\to M$.

\begin{prop}\labell{prop:Bb}
{\rm (i)} For all $k\ge 1$ there is an exact sequence
$$
\pi_k\Bb \stackrel{(\ev_B)_*}\longrightarrow
\pi_kM \stackrel{\Tf_*}\longrightarrow \pi_k(B\Diff \TM).
$$
\NI{\rm (ii)} $(\ev_B)_*\bigl(\pi_1(\Bb)\bigr)\subseteq \Tev_*\bigl(\pi_1(M^{\TM})\bigr)$.
\end{prop}

The following lemma explains some of the elementary properties of the  generalized evaluation subgroup
$\im\, \Tev_*$.   The present formulation of (i) is due to K\c edra. The results here can be be considerably extended; see 
 Gottlieb~\cite{Gott2}.

\begin{lemma}\labell{le:Tev} {\rm (i)} Each 
$\al\in \Tev_*(\pi_1(M^{\TM}))$ 
 acts trivially on the elements in the image of $\pi_*(\TM)$ in $\pi_*(M)$.  In particular, $\al$ acts trivially on $\pi_iM$ for $1\le i< 2n-1$, 
 where $2n: = \dim M$.
 \smallskip
 
 \NI {\rm (ii)}  $\Tev_*(\pi_1(M^{\TM}))$
is trivial if $M$ is a blow up.\smallskip
 
  \NI {\rm (iii)} $\Tev_*(\pi_1(\Bb))$
is trivial  if $\chi(M)\ne0$.\end{lemma}
 
 We begin the proofs with the following lemma.
  
\begin{lemma}\labell{le:ee}  Fix a point $\Tp\in \Si$.  
 Any element $\al\in\pi_1(\ov\Ee)$ can be represented by a 
path $g_t:\Si\to \TM$ such that $g_t(\Tp) = \Tp$ for all $t\in S^1$.
 \end{lemma}
 \begin{proof} Lift $\al$ to a path
 $g_t:\Si\to \TM, t\in I,$ in $\Ee$.  Extend it
 to a path  $\Tg_t$, $t\in I$, in $\Diff\TM$ such that $\Tg_0= id$ and consider the mapping torus $\TP\to S^1$ of $\Tg_1$:
  $$
  \TP: = \TM\times I/\bigl((\Tx,1)\!\sim\!(\Tg_1(\Tx),0)\bigr).
  $$
  Since $\Tg_1$ lies in the identity component of $\Diff\TM$, this is a smoothly trivial bundle and hence its homology is the product
  $H_*(\TM)\otimes H_*(S^1)$.
  
Let $\Si_t: = g_t(\Si)$ and denote by $\TSi\subset \TP$ the union of the submanifolds $\Si_t\times \{t\}, t\in S^1$.  Then 
$Y_t: = \Si_t\cap \Si$ is nonempty for all $t$.  Moreover, by jiggling the $\Tg_t$ we may arrange that 
$\TY: = \cup_t Y_t\times \{t\}\subset \TP$ is a manifold of dimension
$2n-3$ and that $Y_0$ is connected and contains $\Tp$.  We claim that there is a smooth map
$\ga:S^1\to  \TY$ whose projection $\pi\circ\ga$ onto $S^1$ has degree $1$.  
To see this, note first that $\pi^*(dt)$ is nonzero on $\TY$ because
it is Poincar\'e dual (in $\TY$) to $[\TY\cap \pi^{-1}(t)] = 
[\Si]\cap [\Si]\ne 0$.  Hence $\pi^*(dt)$ does not vanish on some cycle $\ga:S^1\to  \TY$ and we can arrange that its projection to $S^1$ has degree $1$ because $Y_0$ is path connected.
Let $\be: = \pi\circ\ga: S^1\to S^1$.  Then
$$
\ga(t) \in Y_{\be(t)},\quad t\in S^1.
$$
Replacing $\Tg_t$ by the homotopic loop $\Tg_{\be(t)}$, we can 
assume that
$\ga(t)\in Y_t\subset \Si$ for all $t\in S^1$. But  $\ga$ contracts in $\Si$ to the point $\Tp$.  Therefore we can homotop the loop $g_t(\Si)$ in $\ov\Ee$ so that 
$\Tp\in g_t(\Si)$ for all $t$,
and then choose $g_t$ to fix $\Tp$. \end{proof}

  \NI
 {\bf Proof of Proposition~\ref{prop:Bb}.}
Let $\TP\to S^k$ be the blow up of $M\times S^k$ along $s_\si$, where $\si:S^k\to M$. 
Suppose that $\pi:\TP\to S^k$ is smoothly trivial.  Then there is a smooth family of diffeomorphisms
$\Tg_t: \TM\to \TM_t: = \pi^{-1}(t)$ for $t\in S^k$ such that $\Tg_0 = id$.  Define 
$\phi_t\in \Bb$ by 
$$
\phi_t: = pr_M\circ\phi_P\circ \Tg_t,\quad t\in S^k,
$$
 where $\phi_P:\TP\to M\times S^k$ is the blow down map and $pr_M$ is the projection to $M$.  By construction $\ev_B(\phi_t) = \si(t)\in M$,
 for all $t\in S^k$. Therefore $\ker \Tf_*\subset \im (\ev_B)_*$.

 To prove the converse, suppose that $[\si] = (\ev_B)_*([\phi_t])$.
 Then 
 $$
 \Phi: \TM\times S^k\to M\times S^k,\quad (x,t)\mapsto (\phi_t(x),t),
 $$
 is a smooth map that is a diffeomorphism over $M\times S^k\setminus s_\si$ and is a blow down over each point in the section $s_\si$.
 Hence $P\cong  \TM\times S^k$ by the uniqueness of the blow up
 construction.
This proves (i).

 Now suppose that $k=1$.
Let $\Si_t'\subset \TM_t, t\in S^1,$ be the copy of the exceptional divisor 
that is blown down by $\phi_P$ and define  $\Si_t: = 
\Tf_t^{-1}(\Si_t')$.
This is a loop in $\ov\Ee$.
Choose a  family of decompositions  
 $\phi_t: = h_t\circ\phi_0\circ \Tg_t, t\in I,$ for the loop 
 $\phi_t\in \Bb$,
 starting with $h_0=\Tg_0=id$.  Then $\Tg_t^{-1}(\Si) = \Si_t$.
 Lemma~\ref{le:ee} implies that
  we may assume that $\Tg_t(\Tp) = \Tp$ for all $t$.
 Therefore $\ev_B(\phi_t): = h_t(p) = \Tg_t(\Tp) =:\Tev(\phi_t)$.  This proves (ii).\QED
 
 \NI
 {\bf Proof of Lemma~\ref{le:Tev}.}  $\al\in \pi_1M$ acts trivially on $\pi_i(M,p)$ if every based map $f:(S^i,x_0) \to 
 (M,p)$ extends to a map $F:
 S^i\times S^1\to M$ such that $F|_{pt\times S^1} = \al$.
But  if $\al = \Tev_*([\phi_t])$  we may define $F: S^i\times S^1\to M$ by setting
$$
F(x,t) = \phi_t(\Tf(x)),
$$
where $\Tf:(S^i,x_0)\to (\TM,\Tp)$ lifts $f$, i.e. satisfies
$\phi_0\circ\Tf = f$. This proves the first statement.
 
 Next observe that when $i<2n$ every element in 
 $\pi_i M$ lifts to $\TM$.
 In fact, since  $\TM$ is diffeomorphic to 
   the connected sum of 
   $M$ with $\C P^n$ (oppositely oriented), $\TM$ has a 
  cell decomposition with
  one vertex at $\Tp$, one $2n$-cell $e^{2n}$, and 
  $(2n-1)$-skeleton 
  of the form $X\vee\Si$, where $X$ is the  $(2n-1)$-skeleton of $M$.
  This proves (i).  

 (ii) holds because when $M$ is the  blow up $\TX$
of $X$ every nonzero element of
$\pi_1M$ acts nontrivially on $\pi_2(M)$. This follows immediately from the fact that $M $ is the connected sum of $X$ with $(\C P^n)^{opp}$, 
so that the universal cover of $M$ is the connected sum of
the universal cover of
$X$ with $|\pi_1(X)|$ copies of  $(\C P^n)^{opp}$ 
on which the group of deck transformations acts effectively.

To prove (iii), consider
$$
Z: = gr\, \phi_0 = \bigl(\Si\times p\bigr)\cup Z_0\subset \TM\times M,
$$
where $Z_0 = \{(\Tq,\phi_0(\Tq)): \Tq\in \TM\setminus\Si\}$.
Let $h:M\to M$ be a  diffeomorphism $C^1$ close to the identity,  with nondegenerate fixed points and such that $h(p)\ne p$.  Then $Z$ intersects the graph of $\phi_0': = h\circ\phi_0$ transversally in $\chi(M)$ points.
(Since $Z$ is a stratified space with top stratum of half the dimension of the ambient manifold and 
one other stratum of codimension $2$, it makes sense to talk of transversality here.  Note that all the intersection points must lie in $Z_0$.)  
Thus
$[Z]\cdot[Z] = \chi(M) \in H_0(\TM\times M)$.  

Now suppose that $\phi_t, t\in S^1,$ is a loop in the space of smooth maps $\Map^\infty(\TM,M)$ based at $\phi_0$.  We may perturb it to a loop $\phi': = \{\phi_t'\}$, based at $\phi_0': = h\circ\phi_0$ (where $h$ is as above) and such that
its graph in $\TM\times M\times S^1$ is transverse to $Z\times S^1$. 
We will show that $\Tev_*(\phi)=\Tev_*(\phi')$ vanishes by finding a 
loop $\Tp_t\in \TM$ of {\it coincidence points}, i.e. points where
\begin{equation}\labell{eq:co}
p_t: = \phi_t'(\Tp_t) = \phi_0(\Tp_t),\qquad t\in S^1.
\end{equation}
(In fact, we will achieve this after a further deformation of $\phi'$.)  Then the 
 loop $t\mapsto \phi_t'(\Tp_t)$
 is homotopic in $(M,p_0)$ to the sum of the loops
 $t\mapsto \phi_0(\Tp_t)$ and $t\mapsto \phi_t'(\Tp_0)$.
 Since equation (\ref{eq:co}) implies that 
 the first two loops are equal, the third loop $t\mapsto \phi_t'(\Tp_0)$ is contractible. 

To find the coincidence points we use smooth intersection theory (which is why we restrict to the space $\Bb$ of smooth maps.)
It is convenient to take $t\in [0,1]=:I$ instead of $t\in S^1$.
Note that for generic $t\in I$, the signed number of points in the 
intersection of $
gr\,\phi_t'$ with $Z\times {t}$ is $\chi(M)\ne 0$.
  Hence, the intersection of $gr\,\phi'$ with
 $Z_0\times I$ contains at least one connected component 
 that 
 projects to $(I,\p I)$ by a map of degree $1$.  Reparametrize $\phi'$ so that this projection is the identity.
 Then, because all intersection points lie in $Z_0\times S^1$, 
there is a path $\Tp_t, t\in I,$ in $M\setminus \{p\}$ such that $\phi_t'(\Tp_t) = \phi_0(\Tp_t)$ for all $t$. 
We now extend the path $\Tp_t, t\in I,$ to a loop $\Tp_t, t\in [0,2]$,
by choosing $\Tp_t, t\in [1,2]$ to be any path in $\TM$ from $\Tp_1$ to $\Tp_0$.
Since $\phi_0'=\phi_1' = h\circ \phi_0$, $\Tp_0$ and $\Tp_1$ are fixed by $h$.  We may assume there is an isotopy $h_t$ from $h_0=id$ to $h_1=h$ that also fixes $\Tp_0,\Tp_1$.  Then extend $\phi'$ by the contractible loop $\phi_t': = h_{\be(t)}\circ \phi_0, t\in [1,2],$ where $\be(1)=\be(2)=1$ and $\be(t) = 0$ for $t\in [1+\eps,2-\eps]$.
After reparametrizing $\Tp_t, t\in [1,2]$,  Eq. (\ref{eq:co}) is satisfied for $t\in [0,2]$.
%%DD Since $\phi_t' = \phi_0$ for $t\ge 1$, we can achieve (\ref{eq:co})
%%DD by a further reparametrization.
%
% and then extend the parameter $t$ over the interval $[1,2]$, setting
% $\phi_t' = \phi_0'$ for $t>1$.
% Because all intersection points lie in $Z_0\times S^1$, 
% we may therefore find a reparametrization map $\be:(S^1,0)\to (S^1,0)$ of degree $1$ and a loop $p_t, t\in S^1,$ in $M\setminus \{p\}$ such that $\phi_t'(\Tp_t) = p_t$ for all $t$, where
% $\Tp_t: = \phi_0^{-1}(p_t)$.  But the loop $t\mapsto \phi_t'(\Tp_t)$
% is homotopic in $(M,p_0')$ to the sum of the loops
% $t\mapsto \phi_0'(\Tp_t)=p_t$ and $t\mapsto \phi_t'(\Tp_0)$. It follows that the loop $t\mapsto \phi_t'(\Tp_0)$ is contractible.  
% The result follows.
\QED 

\begin{rmk}\rm It is not clear how to extend the arguments
in this section  to the case $k>1$.  Even if one could prove some analog of  Proposition~\ref{prop:Bb} (ii)  for $k>1$, the map $\Tev_*: \pi_k(M^\TM)\to \pi_k M$ is not easy to work with.  In particular, 
since $M^\TM$ is not an $H$-space, we cannot immediately 
claim that the image of this map is torsion when $k$ is even (as is the case for the usual evaluation map $M^M\to M$.) All the basic questions can of course be phrased in terms of the space $\Ee/U$, 
but this seems no more accessible. 
In the next section we develop some  
homological tools that work when $M$ is $c$-symplectic and $k=2$.
One of their advantages is that they adapt easily 
to the case of many blowups; see \S\ref{ss:kfold}.  
\end{rmk}  

%%%%%%%%%%%%%%%%%%%%%%%%%%%%%%%%%%%%%%%%%%%%%%%%%%%%%%%%%%%%%%%%%
 \section{The $c$-symplectic case}
%%%%%%%%%%%%%%%%%%%%%%%%%%%%%%%%%%%%%%%%%%%%%%%%%%%%%%%%%%%%%%%%%
  
 This section proves Theorem~\ref{thm:csymp},
 Proposition \ref{prop:diff} and Corollary~\ref{cor:csymp}.  Unless explicit mention is made to the contrary, in this section $(M,a)$ is a connected and simply connected closed $c$-symplectic manifold.

 Denote by $\TQ\to B\Diff^U(M,p)$ the pullback of the universal 
 $\TM$-bundle over $B\Diff\TM$ by the blow up map $B\Diff^U(M,p)\to B\Diff\TM.$  It is obtained by blowing up the universal $M$-bundle over $B\Diff^U(M,p)$ along its canonical section.  
For each  map $\si: S^2\to B\Diff^U(M,p)$, define 
$
(\la,\ell): = \bigl(I_a([\si]), I_c([\si])\bigr).
$
Denote by 
$$
\TP_{\la,\ell}\to S^2
$$
the pullback of $\TQ\to B\Diff^U(M,p)$ by $\si$, and by
 $P_{\la,\ell}\to S^2$ the corresponding $M$-bundle with section $s_\si$.
 
The following result clearly implies 
 Theorem~\ref{thm:csymp}.

  \begin{prop}\labell{prop:S2}  Let $(M,a)$ be a simply connected $c$-symplectic manifold. 
  Then the bundle $\TP_{\la,\ell}$ is smoothly trivial 
   only if $(\la,\ell) = (0,0)$.
     \end{prop}
  \begin{proof} Consider the $M$-bundle $P: = P_{\la,\ell}\stackrel{\pi}\to S^2$, and
denote its  coupling class by $\Ha\in H^2(P)$.  This
is the unique class that restricts to $a$ on the fiber $M$ and
is such that the
fiberwise integral $\pi_!(a^{n+1}): = \int_M {\Ha}^{n+1}\in H^2(S^2)$ vanishes.
Choose an area form $\be$ on $S^2$ so that 
the product class  $\Ha + \pi^*([\be])$ evaluates positively over the section $s_\si\subset P$. 
Then the normalization condition on $\Ha$ implies that
$$
\vol\bigl(P,\Ha + \pi^*([\be])\bigr) = \vol(M,a) \,\int_{S^2}\be\; =:\;V \mu_0,
$$
where $\mu_0:=\int_{S^2}\be$ and $V: = \vol(M,a) = %\displaystyle{
\frac 1{n!}\int_M a^n$.
Further
\begin{equation}\labell{eq:la}
\int_{s_\si} \Ha + \pi^*\be = \mu_0 + I_a([\si]) = \mu_0 + \la.
 \end{equation}

Now choose a representative $\Om$ of the class
$\Ha + \pi^*([\be])$
that restricts in some neighborhood $U$ of $s_\si$ to a symplectic form that induces the given almost structure on its normal bundle.\footnote
{It is not necessary to use symplectic geometry at all in this proof. 
We do so because it is well understood how symplectic forms behave under blow up.}
Next choose $\eps_0$ so that a  neighborhood of $s_\si$ of capacity $\eps_0$ (and hence radius $r_0: = \sqrt{\eps_0/\pi}$) embeds 
symplectically in $U$.  Then, for any $\eps<\eps_0$  we can perform the $\eps$-blow up along $s_\si$ symplectically in $U$; more details are below. 
Denote  the resulting bundle by $(\TP_{\la,\ell},\Tilde\Om_\eps)$.
We show below that
\begin{equation}\labell{eq:vol}
\vol(\TP_{\la,\ell},{\Tilde\Om_\eps}) = \mu_0V -  
v_\eps\left(\mu_0 + \la -\frac\ell{n+1} \eps\right),
\end{equation}
where  $v_\eps: 
= \frac{\eps^n}{n!}$ is the volume of a ball of capacity $\eps$.

Note that the underlying smooth bundle $\TP\to S^2$ does not depend on the choice of $\eps$.  Therefore, if $\TP_{\la,\ell}\to S^2$
 were a smoothly trivial fibration, then, for all $\eps$, the volume of $\TP_{\la,\ell}$ 
would be the product of $V-v_\eps: = \vol(\TM,\Tom_\eps)$ with the \lq\lq size"
 $\mu$ of the base. Since $\mu$ could be measured by  integrating  $\Tilde\Om_\eps$ over a  section of $\TP$ which is the same for all $\eps$, $\mu= \mu_1 + k\eps$ would be a linear function of $\eps$.  Therefore,  the two polynomial functions 
$\vol (\TP_{\la,\ell},{\Tilde\Om_\eps})$ and $(V-v_\eps)(\mu_1 + k\eps)$ 
would have to be equal.
This is possible only if $\la =0$ and also $k=\ell=0$.  
The result follows.

  It remains to derive the formula for $\vol(\TP_{\la,\ell},{\Tilde\Om_\eps})$.
   First assume that $\ell = 0$.  Then the section 
   $s_\si$ has trivial (complex) normal bundle in $M\times S^2$. 
   Hence it has a  neighborhood $U_\eps\subset U$ that is symplectomorphic to a product $B^{2n}(\eps)\times s_\si$.
    Thus, by equation (\ref{eq:la}),
the volume of $U_\eps$ with respect to the form $\Om$ is $\vol(U_\eps) = v_\eps(\mu_0+\la)$.  
We may obtain the blow up by cutting out
$U_\eps$ from  $(P, \Om)$ and identifying the boundary via the Hopf map; see Lerman~\cite{L}.   Then 
$$
\vol(\TP_{\la,0},{\Tilde\Om_\eps} ) = \mu_0V -  v_\eps(\mu_0 + \la),
$$
as claimed.

Now consider the case when $\ell \ne 0$.  
Then  the normal bundle to
$s_\si$ in $M\times S^2$ is isomorphic to the product $\C^{n-1}\oplus L_{\ell}$, where $L_\ell\to S^2$ is the holomorphic line bundle with $c_1=\ell$.  Therefore, we can choose $\Om$ so that it restricts in some neighborhood of $s_\si$ to the product of 
a ball in $\C^{n-1}$ with 
 a $\de$-neighborhood  $\Nn_\de(L_\ell)$ 
of the zero section of $L_\ell$.  Identifying $\Nn_\de(L_\ell)$ with
part of a $4$-dimensional symplectic toric manifold, we can see that its volume is 
$h\de - \ell \de^2/2$, where $h = $ area of zero section and
$\de = \pi r^2$ is the capacity of the disc normal to $s_\si$.
Therefore, since $h = \mu_0+\la$ here, 
\begin{eqnarray*}
\vol(U_\eps) &=& \int_0^{\sqrt{\frac \eps\pi}} \vol(S^{2n-3}(r))\cdot\vol 
(\Nn_{\eps-\pi r^2}(L_\ell))\,dr\\
& = & (\mu_0+\la)\eps^{n}/{n!} - \ell\eps^{n+1}/(n+1)!\\
& = & v_\eps\bigl(\mu_0+\la-\ell \eps/(n+1)\bigr).
\end{eqnarray*}
Everything in the previous calculation remains valid except that we have to add $ v_\eps\ell \eps/(n+1)$ to the volume of $\TP_{\la,\ell}$.  This completes the proof.
  \end{proof}
  
 \NI
{\bf Proof of Corollary~\ref{cor:csymp}.}
The above proof of Theorem~\ref{thm:csymp} uses the fact that $\pi_1(M)=0$
 to assert the existence of a coupling class $\Ha$ on 
 the initial  $M$-bundle $P\to S^2$, from which we construct
 the form $\TOm_\eps$
on the blow up $\TP_{\la,\ell}\to S^2$.  In the situation of this corollary we
are starting
with  a  trivial bundle $M\times S^2\to S^2$ which always has a coupling class, namely $pr_M^*(a)$.  Hence the required form
$\TOm_\eps$ exists on
 the blow up  $\TP_{\la,\ell}\to S^2$. 
The previous argument now applies to show that 
this bundle is nontrivial whenever 
$(\la,\ell)$ is nonzero. The corollary now follows as in the case when $\pi_1(M)=0$.\QED

\NI
{\bf Proof of Proposition~\ref{prop:diff}.}  
Let $X: = \C  P^2\#
  {\ov{\C P}}\,\!^2$, the one point blow up of $\C P^2$.
   Let $a\in H^2(\C P^2)$ be the class of the 
  K\"ahler form $\tau_0$ on $\C P^2$ normalized to integrate to $1$ over the line.
   By \cite{LM2}, 
  any symplectic form $\tau$ on $X$ is
  diffeomorphic to a form  obtained by blow up from $\tau_0$.
  
  Suppose first that $\tau = \tau_\eps$ is the $\eps$-blow up of $\tau_0$.  
By Abreu--McDuff~\cite{AM}, 
 $\pi_1(\Symp(X,\tau)) = \Z$ is  generated by the circle action that rotates the fibers of the
 ruling $X: = \Ppp(L_1\oplus\C)\to \C P^1$ once.  If $X\to P\to S^2$ is the  fibration corresponding to this generator, then one can blow down the exceptional divisors in the fiber to obtain a $\C P^2$-bundle over $S^2$ with section $s$. Note that $I_c(s)=\pm 1$, and we can choose the sign of the generator so that $I_c(s)= 1$.
 
  The blown down bundle is nontrivial, but has order $3$ since $\pi_1(\Symp(\C P^2)) = \pi_1(PU(3)) = \Z/3\Z$; see \cite[Ch~9]{MS}.
  Thus three times the generator of $\pi_1(\Symp(X,\tau)) = \pi_2(B\Symp(X,\tau))$ is the bundle given by blowing up the product $\C P^2\times S^2$ along the graph of the map
 $\si:S^2\to \C P^2$ with image equal to the positively oriented line. Thus it
  corresponds to the bundle $\TP_{\la,\ell}$
 where $(\la,\ell) = (1,3)$. Therefore, the image of 
 $\pi_2(B\Symp(X,\tau))$ in $\pi_2(B\Diff X)$ corresponds to
  the set of fibrations $\TP_{\ell/3,\ell}$, $\ell\in \Z$.
   
   Now consider the form $\tau': = \phi^*(\tau_\eps)$, 
   where $\phi:X\to X$ is
   complex conjugation.  Then $(X,\tau')$
  is  obtained by blowing up the complex conjugate $(\C P^2, \ov \jmath, \ov\tau_0)$ of the standard $\C P^2$. But the first Chern class of the line in $(\C P^2, \ov \jmath,\ov\tau_0)$ is still three times the value of the normalized K\"ahler form on the line.   Hence the image of 
 $\pi_2(B\Symp(X,\tau))$ in $\pi_2(B\Diff X)$ again
 corresponds\footnote
 {
 Observe that the parameter $\la$ in the bundle $\TP_{\la,\ell}$ is determined by the fixed choice of $c$-symplectic class $a$ on $X$.
 Hence  if $\ov\si: S^2\to (\C P^2, \ov\tau_0)$  is an $\ov\tau_0$-symplectic embedding onto the line then $\int \ov\si^*(\Ha) = -1$ and the corresponding $X$-bundle is $\TP_{-1,-3}$.
 }
to the bundles $\TP_{\ell/3,\ell}$, $\ell\in \Z$. 
  
 When appropriately scaled, every symplectic form on $X$ is diffeomorphic 
 either to $\tau$ or $\tau'$ by a diffeomorphism that acts trivially on cohomology and hence fixes both the symplectic class
 and  the first Chern class. 
Since the volume computation is a calculation in 
 the cohomology of $X$, it is also unchanged by such a diffeomorphism.
It follows that 
 if $\tau''$ is an arbitrary symplectic form on $X$, normalized so that its integral over the line has absolute value $1$,
 the image of $\pi_2(B\Symp(X,\tau''))$ in $\pi_2(B\Diff X)$
 corresponds to a set of bundles over $S^2$ with volume forms
equal to those of $\TP_{\ell/3,\ell}$, $\ell\in \Z$.  

Let $\Aa$ be the subgroup of $\pi_2(\Diff X)$ generated by the images of 
 $\pi_2(B\Symp(X,\tau''))$ as $\tau''$ varies over all symplectic forms on $X$.  
Addition in $\pi_2(B\Diff X)$ corresponds to taking the fiber sum of the corresponding bundles.  Hence, because the volume of
$(\TP_{\ell/3,\ell},{\Tilde\Om_\eps})$ depends linearly on $\ell$, 
any bundle $\TP\to S^2$ given by $\al\in \Aa$  has  volume equal to 
that of some bundle $(\TP_{\ell/3,\ell},{\Tilde\Om_\eps})$.
Hence  the bundles $\TP_{0,\ell}, \ell\in \Z,$ do not correspond to elements of $\Aa$.  This completes the proof. \QED

\begin{rmk}\rm  (i)  It would be interesting to find a way to
establish the nontriviality of the
  elements in $\pi_{k}(\Diff\TM)$ for $k>2$  
 constructed by blowing up the product $M\times S^k$ along 
  $pt\times S^k$ with respect to a complex structure
   given by a nonzero element in
  $\pi_k(SO(2n)/U(n))$.   
  For example if $n\ge 4$ one might use a
  nonzero element from $\pi_6(SO(2n)/U(n))$. \MS 
  
  \NI
  (ii) 
 One might also attempt to use a similar construction with
the manifold  $X=S^2\times S^2$, considering this
as the blow up of the orbifold $Y$ 
given by collapsing the antidiagonal
in $S^2\times S^2$ to a point $y_*$.  Since a neighborhood of $y_*$
is diffeomorphic to $\R^4/\!(\pm 1)$, we may construct a bundle $X\to P\to S^2$ from $Y\times S^2$ by blowing up along the constant section $y_*\times S^2$ with respect to a nontrivial family of almost complex structures
$J_z, z\in S^2,$ on $T_{y_*}Y$. 
For example, we could identify $S^2$ with the unit sphere in $\R^3 \equiv (\R e_1)^\perp\subset \R^4$ (where $e_1 = (1,0,0,0)$) and 
define $J_z$ by setting  it equal to the unique element of 
$SO(4)/U(2)$ such that
  $J_z(e_1) = z.$  The resulting loop in $\Diff(S^2\times S^2)$ is certainly not homotopic to a symplectic loop.  However, it seems to represent the sum $[\La] -[\si\La\si]$, 
  where $\si: (z,w)\mapsto (w,z)$
  is the obvious involution of $X$ and $\La$ is the symplectic loop
  that rotates $X$ once about its diagonal and antidiagonal, i.e. if we identify $X$ with the Hirzebruch surface ${\bf P}(L_2\oplus \C)$,
  where $L_2\to S^2$ is the holomorphic line bundle with Chern number $2$, then $\La_t[z_1:z_2] = [e^{2\pi it}z_1:z_2]$. Therefore it is in the subgroup $\Aa\subset \pi_2(\Diff\,X)$ considered above.
  \end{rmk}

%%%%%%%%%%%%%%%%%%%%%%%%%%%%%%%%%%%%%%%%%%%%%%%%%%%%%%%%%%%%%
 \section{The symplectic case}\labell{ss:homol}
%%%%%%%%%%%%%%%%%%%%%%%%%%%%%%%%%%%%%%%%%%%%%%%%%%%%%%%%%%%%%
 
 This section proves Theorem~\ref{thm:homol}.
 We begin with a discussion of characteristic classes defined by 
 parametric Gromov--Witten invariants.  This approach was first suggested by Le--Ono~\cite{LeO} and was then used by Buse~\cite{Bu}.

In the easiest case, we have a class $A\in H_2(M;\Z)$ such that
the formal dimension of the space $\Mm_{g,0}(J,A)$ of unparametrized $J$-holomorphic $A$-curves of genus $g$  is $-k<0$. Let $\langle A\rangle$ denote the orbit of this class under the action of $\pi_0(\Symp\,M)$.
Then there is a corresponding class $c(\langle A\rangle,g)\in H^k(B\Symp\,M;\Q)$
defined as follows.  Because rational bordism equals rational homology, it suffices to evaluate $c(\langle A\rangle,g)$ on cycles of the form 
$\si: Z\to B\Symp\,M$ where $Z$ is a $k$-dimensional closed 
oriented smooth manifold. (For short, we call such cycles smooth.)
 The pullback of the universal $M$-bundle
$M_\Symp\to B\Symp\,M$ 
is a smooth bundle $\pi_\si:P_\si\to Z$.  Then 
$$
\Bigl\langle c(\langle A\rangle,g),\,\,\si_*[Z]\Bigr\rangle: = \PGW^{P_\si,Z}_{g,0},
$$
the \lq\lq number" of unparametrized holomorphic 
$\langle A\rangle$-curves of genus $g$
in the fibers of $P_\si$.  If $\Symp_H$ denotes the subgroup of $\Symp$ that acts trivially on $H_2(M;\Q)$ there is a similar 
class $c(A,g)\in H^k(B\Symp_H\,M)$ that only counts curves in class $A$.\footnote
{
If $Z$ is symplectic and $\pi_1(M)=\pi_1(Z)=0$ then 
Lemma~\ref{le:coup} shows that the fiberwise symplectic forms $\om_z$  extend 
to a symplectic form $\Om$ on $P_\si$. In this case, $H_2(M;\Q)$ injects into $H_2(P_\si;\Q)$ and 
 when $g=0$ the parametric GW invariants are just the usual GW invariants for classes $A\in H_2(M;\Q)\subset H_2(P_\si;\Q)$; cf. \cite[Rem.6.7.8]{MS}.}

Invariants of the above type were used in \cite{LeO,Bu}.
More generally, one might want to count curves that satisfy certain  constraints that are pulled back from the 
cohomology of the universal $M$-bundle $M_{\Symp}\to B\Symp$. 
In the problem at hand we cannot use arbitrary cohomology
classes on $M_{\Symp}$ since we need there to be a relation between the characteristic classes for $M$ and those for $\TM$.  
We therefore consider two different kinds of constraints.

First, we use constraints given by classes in $H^*(M)$ that have canonical extensions to 
$H^*(M_{\Symp})$.  The most obvious classes with this property are the Chern classes  of the  tangent bundle of $M$, which extend to the Chern classes of the  tangent bundle of the fibers of $M_{\Symp}\to B\Symp$.
However, these do not behave well under blowup and so are not very useful here.
%; cf. the proof of Proposition~\ref{prop:homolo2} below. 

Another 
possibility is to use the classes of $\Cc^*$, the subring of $H^*(M)$ generated by $\Cc$.  If $H^1(M;\R)=0$ these classes extend to  
the universal bundle $M_{\Symp_0}\to B\Symp_0 M$ by Lemma~\ref{le:coup}.  For general $M$, they extend to 
 the universal Hamiltonian bundle $M_{\Ham}\to B\Ham\, M$ by \cite{LMP,LM}. The proof of Lemma~\ref{le:coup} shows that in both cases one can pick a unique normalized extension to $M_{\Ham}$.
 For each $\ka\in \Cc^*$, choose a representing polynomial
 $q_{\ka}: = q_{\ka}(a_1,\dots,a_m)$ in the elements $a_i\in \Cc$.
 Then  $\ka$ extends to $\Hq_{\ka}: = 
 q_{\ka}(\Ha_1,\dots,\Ha_\ell)$.  (This extension {\it depends}
  on the choice of polynomial $q_{\ka}$.)
 We then get classes of the form
\begin{equation}\labell{eq:chc}
 c(A,g; \Hq_{\ka_1},\dots,\Hq_{\ka_\ell})\in H^{k+d-2\ell}(B\Ham\,M;\Q),\quad d: = \sum_i \deg \ka_i,
\end{equation}
where $-k= 2n+2c_1(A) -6\ge0$. The value of such a class on $\si_*[Z]$ is (intuitively) the number of $\ell$-pointed genus $g$ $A$-curves in the fibers of $P_\si\to Z$ meeting
 cycles representing the Poincar\'e duals (with respect to $P_\si$) of the classes
 $\Hq_{\ka_1},\dots,\Hq_{\ka_\ell}$. 
 
 Second, we consider the case when $(M,[\om])$ has the $c$-splitting
 property, i.e. the rational cohomology of the
  total space  $M_{\Ham}$ of the universal 
 $M$-bundle over $B\Ham \,M$ is additively isomorphic to 
 $H^*(M)\otimes H^*(B\Ham\, M)$; cf. \cite{LM}. For example,
 by Blanchard \cite{Bl} the hard Lefschetz property mentioned in Theorem~\ref{thm:homol}
 implies the $c$-splitting property.  In this case, all elements of $H^*(M;\R)$ extend to
 $M_{\Ham}$.  Therefore there are classes
\begin{equation}\labell{eq:chc2}
 c(A,g; \Hb_1,\dots,\Hb_\ell)\in H^{k+d-2\ell}(B\Ham\,M;\Q),\quad d: = \sum_i \deg \Hb_i,
\end{equation}
that depend on the chosen extensions $\Hb_i$ of the classes $b_i\in H^*(M;\R)$.

 The proof that these classes are well defined when $g=0$ and
 $(M,\om)$ satisfies a suitable semi-positivity hypothesis is given in Buse~\cite{Bu} and may also be extracted from~\cite{MS}.  The 
 general case is established by the usual methods of the virtual moduli cycle; cf. Li--Tian~\cite{LiT} for example.

 We now return to the specific case at hand. 
  In our examples, $g=0$ and so we will often omit it from the notation.  Let us first suppose that
   $[\si]$ is detected by some element $\ka\in\Cc^*$.   Let $-e\in H^2(\TM)$ be Poincar\'e dual to the exceptional divisor,
 so that $\PD(e^{n-1}) = E,$ the class of a line in the exceptional divisor.  Denote by $\Tka: = \phi_M^*(\ka)$ the pullback of $\ka\in \Cc^*_M$ to $\Cc^*_{\TM}$.
Then %for every  $\Tka\in  (H^2)^k(\TM;\Q)$ 
 the characteristic class of the form (\ref{eq:chc})
 $$
 c(E; \He^{\,n-1}, \He^{\,n-1}\,\Hq_\Tka) \in H^{2k}(B\Symp_0\TM)
 $$
 is well defined.  As we see below, when the bundle $\TP\to Z$ is constructed as a blow up, the value of 
 $ c(E; \He^{\,n-1}, \He^{\,n-1}\,\Hq_\Tka)$ on $\si_*[Z]$
 is given by integrating
 $\Hq_\Tka$ over a suitable section of the bundle $\TSi\to Z$ of 
 exceptional divisors.  
 
 The basic reason why the 
 next proposition holds is that the extensions $\Hq_\ka$ are well behaved under blow up. More precisely, suppose that the $\TM$-bundle $\TP\to Z$ is the blow up of 
the trivial bundle
$P=M\times Z\to Z$ along $gr_\si$, and let $\phi_P:\TP\to P$ be the blow down.  Then, if $H^1(M;\Q)=0$, it is easy to see that
 $\phi_P$ pulls the coupling class $\Ha \in H^2(P)$ of $a\in \Cc_M$ back to the coupling class of $\phi_M^*(a)$ in $\TP\to Z$. (One just needs to check this when $\dim Z=2$ and then the normalization condition for the coupling class $\Hc$ is that $\Hc^{n+1} = 0$.)  
 Hence 
\begin{equation}\labell{eq:pb}
 \phi_P^*(\Hq_\ka|_P) = \Hq_{\Tka}|_{\TP}\quad\mbox{where } \Tka: = \phi_M^*(\ka).
\end{equation}
 
 \begin{prop}\labell{prop:homolo} Assume that  $H^1(M;\R)=0$. Let $\si:Z\to M$ be a smooth cycle representing the
 homology class  $[\si]\in H_{2k}(M)$, and denote by 
 $\Tf_\eps:Z\to B\Symp_0 (\TM,\Tom_\eps)$ the cycle given by the
 $\eps$-blow up of
 $M\times Z\to Z$ along $gr_\si$.
   Then, for  every polynomial 
   representative $q_\ka$ of $\ka\in \Cc^*$,
 $$
\Bigl\langle  c(E; \He^{\,n-1}, \He^{\,n-1}\,\Hq_\Tka),
\,\,(\Tf_\eps)_*[Z]\Bigr\rangle  =\ka([\si]).
 $$
 \end{prop}
 
 The proof uses the following lemma.
 
 \begin{lemma}\labell{le:GW0} Let $E\in H_2(\TM;\Z)$ be the class of a line in the exceptional divisor. Then 
 $\GW^{\TM,E}_{0,2}(e^{n-1},e^{n-1}) = 1$.
\end{lemma}
\begin{proof}    
Choose $J_0$ on $\TM$ that 
equals the standard integrable structure near  the exceptional divisor $\Si$.
Let $\ell_1,\ell_2$ be two disjoint
  embedded $2$-spheres in $\TM$ representing $E$ that meet $\Si$ once transversally (with sign $-1$). Consider the moduli space $\Mm_{0,2}(E,J_0)$ 
  of $J_0$-holomorphic $E$-spheres
   with $2$-marked points in $\TM$ (quotiented out by $\C$, the reparametrizations that fix $0,\infty\in S^2$.) 
   Since the exceptional divisor $\Si$ is $J_0$-holomorphic and $E\cdot \Si=-1$, the only $J_0$-holomorphic $E$-curves in $\TM$ are the obvious curves in $\Si$.
    Thus $\Mm_{0,2}(E,J_0)$ is compact, except that the two marked points can come together.  Hence it can be completed to a {\it manifold} by adding in the relevant stable maps, that have a ghost bubble with the two marked points on it. We  call the compactified space $\oMm_{0,2}(E,J_0)$.

  All the elements in $\oMm_{0,2}(E,J_0)$ are regular by construction.
  So $\oMm_{0,2}(E,J_0)$ is a manifold of dimension
  $2n + 2c_1(E) - 2 = 4n-4$.  Consider the evaluation map
  $$
  \ev: \oMm_{0,2}(E,J_0)\to \TM\times \TM.
  $$
  It meets $\ell_1\times \ell_2$ exactly once and transversally.
  This proves the lemma. \end{proof}
  
  \NI{\bf Proof of Proposition~\ref{prop:homolo}.}
  Consider the bundle $\pi_Z:(\TP, \TOm_\eps)\to Z$ given 
  by the $\eps$-blow up of
 $M\times Z\to Z$ along $gr_\si$.
  Choose a family $J$ of fiberwise $\TOm_\eps$-tame 
  almost complex structures $J_z$, $z\in Z$, on 
  $\TP$, that are standard near each exceptional 
  divisor $\Si_z$. Denote by $\TSi$ the union of 
  the exceptional divisors.  Then the moduli space $\Mm_{0,0}(E,J)$ of unparametrized $J$-holomorphic $E$-spheres is a compact manifold
  of dimension $4n-4+2k$ consisting of regular curves.
The evaluation map $\ev:\Mm_{0,2}(E,J)\to \TP\times \TP$ is a 
pseudocycle of dimension $4n+2k$.\footnote{Note that we do not 
compactify here, i.e. the $2$ marked points are assumed 
distinct, but we do divide out by the reparametrization group. 
Pseudocycles are discussed in~\cite{MS,Z}.}
  Further   
 $ \bigl\langle  c(E; \He^{\,n-1}, \He^{\,n-1}\,\Hq_\Tka),\,\,
 (\Tf_\eps)_*[Z]\bigr\rangle$ is the intersection number of $\ev:\Mm_{0,2}(E,J)\to \TP\times \TP$ with the product $X_1\times (X_2\cap X_3)$,
 where $X_1, X_2$   represent the Poincar\'e dual of
$\He^{\,n-1}$ and $X_3$ represents that of $\Hq_\Tka$. 
 
  To calculate this, choose  
 a smoothly embedded\footnote
 {
 The Poincar\'e dual of
  $\He$ can be represented by an embedded  submanifold $Q$:
  take the inverse image of a hyperplane by a map $P\to \C P^N$ 
  representing $\He$. Hence we can take $X$ to be an iterated 
  intersection.  Since $[X]\cap [\TM] = E$, roughly speaking 
  $X$ intersects a generic fiber of $\TP\to Z$ in a line.}
   submanifold 
  $X\subset \TP$ that represents the Poincar\'e dual of
  $\He^{\,n-1}$, choose two copies  $X_1, X_2$ of $X$ that intersect
  $\TSi$ in transversally intersecting $2k$-dimensional submanifolds  $Y_1, Y_2$, 
  and consider the submanifold $\Mm^{cut}: = \ev^{-1}(X_1\times X_2)$
  of  $\Mm_{0,2}(E,J)$. 
    Evaluation at the second marked point $\ev_2: \Mm^{cut}\to \TP$ 
   is a $2k$-dimensional pseudocycle (with image in $\TSi\cap X_2$)
   and, by definition, 
   $$
    \bigl\langle  c(E; \He^{\,n-1}, \He^{\,n-1}
    \,\Hq_\Tka),\,\,(\Tf_\eps)_*[Z]\bigr\rangle = 
   \bigl\langle \ev_2^*(\Hq_\Tka),\,\,\Mm^{cut}\bigr\rangle =
   \bigl\langle \ev_2^*\phi_P^*(\Hq_\ka),\,\,\Mm^{cut}\bigr\rangle,
   $$  
   where the second equality holds by (\ref{eq:pb}).

We now claim   that $\phi_P\circ\ev_2:\Mm^{cut}\to P=M\times Z$ represents $[gr_\si]$.  To see this, consider the codimension $2$ submanifold 
 $Z_1: = \pi_Z(Y_1\cap Y_2)\subset Z$.
Because there is a unique $J$-holomorphic $E$-curve through 
  each pair of distinct points in $\Si_z$,
  $$
  \Mm^{cut}\cap \ev^{-1}\bigl(\pi_Z^{-1}(Z\setminus Z_1)\bigr)
  \cong Y_1\times_{Z\setminus Z_1} Y_2.
  $$
  But by construction
 the map $\phi_P:Y_i\to gr_\si$ has degree $1$. Hence 
 $\phi_P\circ\ev_2:\Mm^{cut}\to gr_\si$ also has degree $1$.
 
Therefore,
$$
   \bigl\langle \ev_2^*\phi_P^*(\Hq_\ka),\,\,\Mm^{cut}\bigr\rangle =
      \bigl\langle \Hq_\ka,\,\,[gr_\si]\bigr\rangle = 
      \bigl\langle \ka,\,\,[\si]\bigr\rangle,
$$
where the second equality holds because, for the trivial bundle $P=M\times Z\to Z$, we have $\Hq_\ka= pr_M^*(\ka)$.  \QED\MS
 
 \begin{cor} Theorem~\ref{thm:homol} holds.
\end{cor}
\begin{proof}  This is an immediate consequence of Proposition~\ref{prop:homolo}.
\end{proof} 
 
 We next turn to the proof of Proposition~\ref{prop:homol2}, which
 concerns classes $[\si]$
that can be represented by smooth, simply connected cycles, i.e. smooth maps $\si:Z\to M$ where the manifold $Z$ is simply connected.
The next lemma is no doubt well known. We include a proof for completeness.

\begin{lemma}\labell{le:sc} If $\dim M>4$ then the classes 
represented by smooth, simply connected cycles generate the image
in $H_*(M;\Q)$ of the homology of the universal cover of $M$.
\end{lemma}
\begin{proof}
It is obvious that these classes lie in this image since any such map $\si$ factors through the universal cover.  To prove the converse, it suffices to prove that if $\pi_1(M)=0$ and  $Z\subset M$ is a smooth oriented embedded\footnote
{We can assume that $Z$ is embedded by Thom's well known result that 
rational homology is generated by embedded smooth cycles.} 
$d$-manifold with $d\ge 3$ then one can do surgery on $Z$ 
 to make it simply connected also.   Since $\dim M>4$ 
 every embedded loop $\ga$ in $Z$ is the boundary
  of an embedded $2$-disc $D$ in $M$ that 
  %%D
  we can chose to be nowhere tangent to $Z$ along
  $\p D = \ga$.  The normal bundle $\nu_D$ of $D$ in $M$ splits 
along the boundary as a sum $\eps^{d-1}\oplus \eps^{2n-d-1}$ where
$\eps^{d-1} $ is the trivial normal bundle to $\ga$ in $Z$. If this splitting extends over $D$ so that $\nu_D$ is a sum of oriented bundles
$\xi^{d-1} \oplus \xi^{2n-d-1}$ then we can kill $\ga$ in $\pi_1(Z)$ by doing (possibly twisted) surgery to $Z$, removing a neighborhood of $\ga$ and adding the boundary of a neighborhood of $D$ in $\xi^{d-1}$.
%%D
Therefore we need to check that
the element $\al$ of $\pi_2(B\SO(2n-2))$ represented by
the normal bundle to $D$ with its given boundary trivialization
lifts to an element of $\pi_2(B\SO(d-1)\times B\SO(2n-d-1))$.  The obstruction is the image of $\al\in \pi_1(Gr)$, where $Gr: = 
\SO(2n-2)/\SO(d-1)\times \SO(2n-d-1)$ is the Grassmannian of oriented $(d-1)$-planes in $\R^{2n-2}$. 
%Therefore we have a lift $\p D\to B\SO(d-1)\times B\SO(2n-d-1)$
%of the classifying map $(D,\p D) \to (B\SO(2n-2),*)$ of the normal bundle (with its given boundary trivialization) 
%and need to extend this over $D$. 
 Since $\pi_1(Gr) = 0$,
the extension is unobstructed.
  \end{proof}
 
We now restate Proposition~\ref{prop:homol2} for the convenience of the reader.
  
 \begin{prop}\labell{prop:homolo2} 
 Let $\si:Z\to M$ be a smooth and simply connected nontrivial cycle,
 and let  $\Tf_\eps:Z\to B\Ham (\TM,\Tom_\eps)$ be the corresponding classifying map.
 Then $(\Tf_\eps)_*[\si]\ne0$ if one of the following conditions hold.
 \smallskip
 
 \NI{\rm (i)} There is $\ka\in \Cc^*$ such that $\ka([\si])\ne 0$.
   
 \smallskip
 
 \NI{\rm (ii)}  
$Z=S^{2k}$ and $(M,\om)$ has the hard Lefschetz property.
%\smallskip

%\NI
%{\rm (iii)} $Z=S^{2k}$ and the Chern class $c_k$ does not vanish on $[\si]$. 
\end{prop}
 \begin{proof}
 Case (i) has the same proof as Proposition~\ref{prop:homolo}.
 The crucial point is that an analog of Lemma~\ref{le:coup} holds in this setting.  We know that all classes $a\in \Cc$ extend to 
 $M_{\Ham}$ by Lalonde--McDuff--Polterovich~\cite{LMP,LM}.  Moreover, 
 because $\pi_1(Z)=0$  any two extensions differ by a class pulled back from the base.  Hence the normalization condition is enough to provide a {\it unique} extension that satisfies the compatibility condition~(\ref{eq:pb}).
 
 In case (ii) all classes in $H^*(M;\Q)$ extend to $M_{\Ham}$ by 
 Blanchard \cite{Bl}.  The same holds for $\TM$ since this also has the hard Lefschetz property.  The difficulty is to ensure that
 the extensions satisfy the analog of (\ref{eq:pb}). We do this by imitating the normalization procedure in Lemma~\ref{le:coup}.
 Given $b\in H^d(M)$ we claim there is
$\Hb\in M_{\Ham}$ such that
\begin{equation}\labell{eq:norm}
\pi_!(\Hb\,\Ha^n): = \int_M\Hb\,\Ha^n = 0\in H^d(B\Ham),
\end{equation}
where $\Ha$ is the coupling class (i.e. normalized extension of
the symplectic class $[\om]$), and $\int_M$ denotes integration over the fiber.  (The proof is as in Lemma~\ref{le:coup}; given any extension $\Hb\,'$ simply subtract from it the pullback of a  suitable
multiple of $\pi_!(\Hb\,'\,\Ha^n)$.)  This normalized extension $\Hb$ 
may not be unique.  However, its pullback to any bundle $P\to Z$ 
over a sphere $Z=S^d$ is unique.  Moreover because the normalization condition on $P$ is simply that $\Hb\,\Ha^n=0 \in H^{2n+d}(P)$
the analog of~(\ref{eq:pb})
is also satisfied. Thus, if $Z=S^d$ and $b\in H^*(M)$,
$$
\phi_P^*(pr_M^*(b)) = \Hb|_{\TP},
$$
where $\Hb|_\TP$ denotes the restriction to $\TP$ of {\it any} extension of $\phi_M^*(b)$ to $\TM_{\Ham}$ that satisfies 
(\ref{eq:norm}). 

Now argue as in  
Proposition~\ref{prop:homolo} using the characteristic class
$c(E; \He^{\,n-1}, \He^{\,n-1} \,\Hb)$ instead of
$
c(E; \He^{\,n-1}, \He^{\,n-1}\,\Hq_\Tka)$,  where $b$ is chosen so that $b([\si])\ne 0$ and $\Hb$ is any extension of $\phi_M^*(b)$ that satisfies (\ref{eq:norm}).
This proves  (ii).
\end{proof}

\begin{example}\labell{ex:tor} \rm
 Consider $M=\T^{2n}$ with a standard form. 
Because $\ev:\pi_1(\Symp\, M)\to \pi_1(M)$ is surjective in this case, 
the bundle $\TM\to \TP\to M$  of equation (\ref{eq:MM})
is symplectically trivial over the $1$-skeleton of the base $M$.  Therefore, because $\pi_1(M)$ has no torsion, this bundle is classified by a map into 
$B\Ham\,\TM$; see
\cite[Erratum]{LM}.  Because $\T^{2n}$ has the hard Lefschetz property, one might then expect
$H_d(M)$ to inject into $H_*(B\Ham\,\TM)$ for all $1<d\le 2n$.
But this  map is trivial because there is an automorphism of the
product bundle $M\times M\to M$ that takes the diagonal to the constant section; in fact the
 evaluation map
$\ev:\Symp_0(\T^{2n})\to \T^{2n}$ is surjective on homology.
This example explains why one needs some conditions in part (ii)
of  Proposition~\ref{prop:homolo2}, though one might well be able to weaken the ones given.
\end{example}

\begin{rmk}\rm Part (ii) of Proposition~\ref{prop:homolo2} holds whenever $(\TM,\Tom_\eps)$ has the $c$-splitting property for all $\eps<\eps_0$, i.e. whenever the cohomology
Leray--Serre spectral sequence of the universal $\TM$ bundle $\TM_{\Ham}\to B\Ham(\TM,\Tom_\eps)$
 degenerates at the $E_2$ term. We claim that in this case 
$(M,\om)$ also   
has the $c$-splitting property. However,
it is not clear whether the $c$-splitting property
for $(M,\om)$ implies that for $(\TM,\Tom_\eps)$ since there may be symplectomorphisms of $\TM$ that have little to do with those of $M$.

To prove the claim, note that if
 $(\TM,\Tom_\eps)$ has the $c$-splitting property then, because compact subsets of $B\Ham(M,p)$ map compatibly to $B\Ham(\TM,\Tom_\eps)$,
the universal $M$-bundle with section $M_{\Ham(M,p)} \to B\Ham(M,p)$
also has degenerate spectral sequence.  To see that 
the spectral sequence of $M_{\Ham}\to B\Ham(M,\om)$ also degenerates, it suffices to check that the map 
$\pi: M_{\Ham}\to B\Ham(M,\om)$ on base spaces induces an injection on 
cohomology, since then the spectral sequences inject; cf. 
\cite[Lemma~4.1]{LM}.  But this holds because there is a class $\Hc\in 
H^{2n}(M_{\Ham})$ (namely $\Hc: = \Ha^n$ where $\Ha$ is the canonical 
extension of $a: = [\om]$) that restricts to a generator of $H^{2n}(M;\R)$.  Since in this case  the map
$$
H^*(B\Ham;\R)\to H^*(B\Ham;\R),\quad 
b\mapsto\pi_!(\Hc\cup \pi^*(b))
$$
is multiplication by a nonzero constant, $\pi^*$ is injective.
\end{rmk}

%%%%%%%%%%%%%%%%%%%%%%%%%%%%%%%%%%%%%%%%%%%%%%%%%%%%%%%%%%%%%
 \section{Blowing up at many points}\labell{ss:kfold}
%%%%%%%%%%%%%%%%%%%%%%%%%%%%%%%%%%%%%%%%%%%%%%%%%%%%%%%%%%%%%
  
Now suppose that we  blow up more than once. The first result concerns
 simultaneous blow ups.  

  \begin{prop}\labell{prop:csympk}
  Let $(M,a)$ be a $c$-symplectic manifold  and denote by $\TM_k$ its $k$-fold blow  up.  Then:\smallskip
  
  \NI
 {\rm  (i)}  There is a homomorphism $\Tf_*^k: 
  \bigl(\Z \oplus \pi_2M\bigr)^k \to \pi_2\bigl(B\Diff\TM_k\bigr)$ whose kernel is 
contained in the torsion subgroup of $\bigl(\pi_2M\bigr)^k$.\smallskip

\NI
{\rm (ii)}  If $(M,\om)$ is symplectic then there is $\eps_0>0$ such that
the elements in $\Tf_*^k\bigl((\pi_2M)^k\bigr)$ can all be realised in
$B\Symp(\TM_k, \om_\eps)$ whenever the blow up parameter $\eps = (\eps_1,\dots,\eps_k)$ satisfies  $\eps_i\le \eps_0$ for all $i$.
\end{prop}

\begin{proof} Fix a metric on $M$.  Given a (small) constant $\nu>0$,
define  $\De_\nu\subset M^{k}$ by setting
$$
\De_\nu: = \{(x_1,\dots,x_k)\in M^{k}: d(x_j,x_i)\le \nu, i\ne j\}.
$$
Then the product bundle $M\times (M^{k}\setminus \De_\nu)\to  (M^{k}\setminus \De_\nu)$ has $k$  symplectic sections of the form
$(x_1,\dots,x_k)\mapsto (x_j,x_1,\dots,x_k)$ that are 
mutually separated by a distance of at least $\nu$.  Hence, for sufficiently small $\eps_0 = \eps_0(\nu)>0,$
it is possible to blow up $M\times (M^{k}\setminus \De_\nu)$ along these sections simultaneously with any weights $\eps_i\le \eps_0$.
This defines a model  fibration $(\TM_k, \om_\eps)\to (\TP_k, \Om_\eps)\to M^{k}\setminus \De_\nu$ that plays the role of the 
fibration (\ref{eq:MM}).

To prove (ii) 
fix $k$ distinct points $y_j$ in $M$ and
choose smooth maps $\si_\al^j:S^2\to M$,
for $1\le \al\le r$, that represent  a basis for $H_2(M)/{\rm torsion}$ and are 
such that for all $\al, j$ the image of  $\si_\al^j$ does not contain the points $y_i, i\ne j$.
Now set 
 $$
 \nu: = 
\frac 12  \inf\{
d(\si_\al^j(z),y_i):i\ne j, z\in S^2, 1\le \al\le r\}.
 $$
For each $j,\al$ as above, consider the map
$$
\Hat{\si}_\al^j: S^2\to M^{k}\setminus \De_\nu,\quad z\mapsto
(y_1,\dots, y_{j-1},\si_\al^j(z), y_{j+1},\dots y_k).
$$ 
  The element $\Tf_*^k([\Hat{\si}_\al^j])\in \pi_2(B\Symp(\TM_k,\om_\eps))$ is represented by the pullback of the model fibration by $\Hat{\si}_\al^j$.  This defines $\Tf_*^k$ 
  on  a set of generators of $H_2(M^k)/{\rm torsion}$.  To see that it is injective, one calculates volumes as before, 
 considering these as a function of the $k$ blow up parameters $\eps_1,\dots, \eps_k$.  This proves (ii).
 
 The proof of (i) is a straightforward adaptation of that of Corollary~\ref{cor:csymp}, and is left to the reader.
\end{proof}

The above result is weaker than it need be.  Suppose for example that
$H_2(M)$ has rank $r$ and that $k=2$.  Then the argument above detects
 a subgroup of rank $2r+2$ inside $\pi_2\bigl(B\Diff\TM_2\bigr)$.
 However this group should contain a subgroup of rank $2r+4$,
 where the two extra generators $\ga_1,\ga_2$ 
 come from doing one of the blow ups along a line in the exceptional 
 divisor formed by the other blow up. These blow ups can no longer be done simultaneously. Rather, the corresponding fibration $\TP\to S^2$ is a fiber sum $\TP_1\#\TP_2$ where each $\TP_i$ is the 
 $\eps_i$-blow up of $\TM_1\times S^2$ along a suitable section.
 Thus when forming $\TP_1$ one starts with a trivial fibration 
 whose fiber has volume $V-v_{\eps_2}$, while 
 for $\TP_2$ the fiber has volume $V-v_{\eps_1}$.
  More generally, if we blow up $k$ times, then 
  $\TP$ decomposes as $\TP_1\#\cdots\#\TP_k$, where $\TP_i$ is the $\eps_i$-blow up of a trivial bundle with fiber of volume 
  $V-\sum_{j\ne i} v_{\eps_j}$ along a section $s_i$ whose homology class
 lies in a group of rank $r+k-1$ and has one parameter $\ell_i$ describing its normal bundle.   Therefore,
\begin{equation}\labell{eq:volk}
\vol(\TP_{\la,\ell},{\Tilde\Om_\eps}) = 
\sum_i \Bigl(\bigl(V -\sum_{j\ne i} v_{\eps_j}\bigr)\mu_i - 
v_{\eps_i}\bigl(\mu_i + \la_i -\frac{\ell_i}{n+1} \eps_i\bigr)\Bigr)
\end{equation}
where $\la_i: = I_a(s_i)$ is no longer constant but is a sum $\la_i' + \la_i''$, 
where $\la_i' = \sum_{j\ne i} m_{ij}\eps_j$ for $m_{ij}\in\Z$, and $\la_i''\in\R$.
 If the bundle were trivial this polynomial would have the form $(V-\sum_jv_{\eps_j})(\mu +\sum k_j\eps_j)$.
Since $V\ne 0$ and $\vol(\TP_{\la,\ell},{\Tilde\Om_\eps})$ has no terms that are linear in the $\eps_i$ it follows as before that 
all the coefficients $k_i, m_{ij}, \la_i,\ell_i,$ must vanish 
and that $\mu=\sum \mu_i$.  This proves the following result.

\begin{prop}\labell{prop:csympk2}
  Let $(M,a)$ be a $c$-symplectic manifold  
  and denote by $\TM_k$ its $k$-fold blow  up.  Then
  the rank of $\pi_2\bigl(B\Diff\TM_k\bigr)$ is at least 
  $k(r+k)$, where $r: = {\rm rank\,} \pi_2(M)$.
\end{prop}

We now address the question of when these new elements in 
$\pi_2\bigl(B\Diff\TM_k\bigr)$ can be constructed to be symplectic.
 Note that if $(M,\om)$ is a blow up, the size of its exceptional divisor does not constrain the maximal size of a single symplectically embedded ball; for example if one blows up $\C P^2$ just a little, one can 
 blow it up again using a much larger ball.  However the next argument shows that it does constrain the size of suitable families of embedded balls.

\begin{prop}\labell{prop:size}  Let $(X,\om)$ be a symplectic manifold that is itself a blow up of size $\al$, that is $\al$ is the integral of $\om$ over a line in the exceptional divisor $\Si_X$.   
 Suppose that $\si:S^2\to X$ represents a class such that $[\si]\cdot [\Si_X] \in H_0(X)$ is nonzero and construct    
$\pi: (\TP, \Tilde\Om_\eps)\to S^2$  by blowing
up $X\times S^2$ along $gr_\si$. Then:\smallskip

\NI
{\rm (i)}   $\eps < \al$. \smallskip

\NI
{\rm (ii)} If 
$\Tilde\rho_t, 0\le t\le 1,$ is any family of symplectic forms on $\TP$ that start at some blow up form $\Tilde\rho_0: = \Tilde\Om_\eps$ then $\int_{E_\si}\Tilde\rho_1 < \al$,
where $E_\si$ is the class of a line in
the exceptional divisor in the blow up $\TX$.
\end{prop}
\begin{proof} Let $P: = X\times S^2$. We first claim
$$
\GW^{P,E_X}_{0,2}([gr_\si], E_X\times [S^2]) = - [\si]\cdot [\Si_X]\ne 0
$$
where $E_X$ is the class of   a line in the exceptional divisor $\Si_X$.
To prove this, choose an $\om$-tame almost complex 
structure $J$ on $X$ that is integrable near $\Si_X$, Choose an embedded 
$2$-sphere $\ell_X$ in the class $E_X$ that meets $\Si_X$ once transversally, 
and perturb $\si$ so that $gr_\si$, $\Si_X\times S^2$ and  
$\ell_X\times S^2$ are in general position. 
 Then, to each intersection point $x\in gr_\si\cap (\Si_X\times S^2)$ there is precisely one $(J\times j)$-holomorphic 
curve in $X\times S^2$  that meets $\ell_X\times [S^2]$.  
It is regular and,
because $\ell_X\cap \Si_X = -1$, it contributes to the Gromov--Witten invariant with sign opposite to that of the orientation of the intersection at $x$.
This proves the claim.

Now consider $(\TP,\Tilde\Om_\eps)$.  We claim that 
$$
\GW^{P,E_X}_{0,2}([gr_\si], E_X\times [S^2]) = 
\GW^{\TP, E_X-E_\si}_{0,2}(\TSi,E_X\times [S^2]),
$$  
where $\TSi$ denotes the subbundle of $\TP$ formed by the 
exceptional divisors.\footnote
{
This GW invariant (which can be considered as a parametric GW invariant on the bundle $\TP$)
need not come from a characteristic class on $B\Symp_0$ since it is not clear that the submanifold $\TSi$ extends to a codimension $2$ homology class in $\TM_{\Symp_0}$.  For instance,
$\TSi$  is not the Poincar\'e dual of the normalized extension $\He$.}
  This is proved by direct construction.  
Denote by $\Om$ a product symplectic form $\om +\pi^*(\be)$ on $P$, and
 choose an $\Om$-tame almost complex structure $J_P$ on $P$ so that $gr_\si$ is $J_P$-holomorphic and so that $J_P$ is integrable near $gr_\si$.  If $\eps$ is sufficiently small we can perform the blow up of $P$ along $gr_\si$ 
 in such a way that $J_P$ lifts to an $\Tilde\Om_\eps$-tame almost complex structure $\TJ_P$ on $\TP$, i.e. so that the blow down map is $(\TJ_P,J_P)$-holomorphic.  It then follows that there is a bijective correspondence between the curves counted by each of the above Gromov--Witten invariants.
 (For details of a similar proof see \cite[Ch~9.3]{MS}.)
 
 It follows that the integral of $\Tilde\Om_\eps$ over $E_X-E_\si$ must be positive.  This proves (i).  (ii) follows because the Gromov--Witten invariants are unchanged by a deformation of the symplectic structure.
\end{proof}

Now suppose that $B_1, B_2$ are two disjoint symplectically 
embedded balls in $(M,\om)$. Denote by $X_i, i=1,2,$  the blow up of $M$ by $B_i$, and by $X_{12}$ the blow up by both balls.   Denote by $E_i$ the  class of a line in the exceptional divisor
corresponding to $B_i$.  For $i=1,2$ let $\La_i\in\pi_2(B\Diff(X_{12}))$ be the element formed by blowing up $X_i\times S^2$ along the graph of $E_i$.
We would like to be able to conclude from the above proposition
that the element $\La_1 + \La_2$, for example, is not homotopic to an element of $\pi_2\bigr(B\Symp(X_{12},\Tilde\om)\bigr)$  for any blow up form $\Tilde\om$ on $X_{12}$ because we need $\Tilde\om(E_1)< \Tilde\om(E_2)$ to represent $\La_2$ and the reverse inequality to represent $\La_1$.\footnote
{
By \lq\lq blow up form" on a $k$-fold blow up $\TM_k$
we mean a symplectic form on $\TM_k$ that is constructed from $(M,\om)$ by cutting out $k$ disjoint standard balls. Similarly, a symplectic form 
%%D
on the $\TM$-bundle $\TP$ is called a blow up form if it is obtained from a  form on an $M$-bundle over $S^2$ by cutting out a suitable neighborhood of a section.
This is equivalent to requiring that $\TP$ has 
a subbundle with  the exceptional divisor as fiber.
One feature of such forms is that they can always be deformed through symplectic forms so as to decrease the size of the balls.
}
The problem is that Proposition~\ref{prop:size} applies only to bundles $(\TP,\Tilde\Om)$ such that $\Tilde\Om$ can be deformed into a blow up form, while the elements of $\pi_2(\Ham \TM_k)$ correspond to 
bundles  $(\TP,\Tilde\Om)$ where $\Tilde\Om$ is 
an arbitrary extension of a blow up form 
%%D 
on $\TM_k$ 
that may not have such a deformation.  In the next section we show 
that this does not happen in dimension $4$.

%%%%%%%%%%%%%%%%%%%%%%%%%%%%%%%%%%%%%%%%%%%%%%%%%%%%%%%%%%%%%
 \section{The $4$-dimensional case.}\labell{ss:4}
%%%%%%%%%%%%%%%%%%%%%%%%%%%%%%%%%%%%%%%%%%%%%%%%%%%%%%%%%%%%%

This section discusses $k$-fold blow ups of a 
$4$-dimensional manifold $M$.
 First a preliminary lemma. Though well known, 
we include a brief proof for completeness.  (See \cite{MS} for 
further details.)

\begin{lemma}\labell{le:e}  Let $(M,\om)$ be a symplectic $4$-manifold and $E$ the class of a symplectically embedded sphere of self intersection $-1$.  Then if $J_z, z\in Z$, is any generic compact 
$2$-parameter family of $\om$-tame almost complex structures, 
there is a finite subset $\{z_1,\dots,z_k\}\subset Z$
such that $E$ is represented  by a unique
 embedded $J_z$-holomorphic sphere for all $z\in Z\setminus\{z_1,\dots,z_k\}$.  Moreover, 
 for $z=z_i$ the class
$E$ is represented by a nodal $J_z$-holomorphic 
curve consisting of two embedded
spheres, intersecting once transversally and
 with self-intersections $-2$ and  $-1$.
\end{lemma}
\begin{proof}  Because $GW^{M,E}_{0,0}= 1$, $E$ always has a $J$-holomorphic representative.  We must investigate the possible nodal curves representing $E$.  These are finite unions $C_0\cup\cdots\cup C_k$ of spheres in
classes $A_0,\dots, A_k$.
The index $\ind\, D_u$ of a sphere $u:S^2\to M$ in class $A$ in a $4$-manifold
is $4 + 2c_1(A)$, where $D_u$ denotes the linearized Cauchy--Riemann operator.  In a generic $2$-parameter family 
the cokernel of $D_u$ has dimension at most $2$. Hence
if $u$ is simple so that there is a $6$-parameter reparametrization group
this index must be at least $4$, i.e. we must have $c_1(A)\ge 0$. 
 If
$u$ is multiply covered then  the reparametrization group has dimension at least $10$ which means that the index must be at least $8$. In other words $c_1(A)\ge 2$.  Since $\sum_ic_1(A_i) = 1$, it follows that all components
are simple and that all but one (say $A_0$) have $c_1(A_i) = 0$. 
Since 
the components are simple, they satisfy the adjunction inequality:
$c_1(A_i) = 2 + (A_i)^2 - 2d_i$ where $d_i\ge 0$ is 
the defect.\footnote
{
If $C_i$ is immersed, $d_i$ is the number of double points;  $d_i=0$ iff $C_i$ is embedded.}  Hence $(A_i)^2\ge -2$ for $i>0$ and $(A_0)^2\ge -1$.  
Since the nodal curve is connected, there are at least $k$ pairs of components $C_i, C_j$ that  intersect.  Hence, by positivity of intersections, 
 $A_i\cdot A_j \ge 1$ for $i\ne j$, and we find
$$
-1=(A_0+\cdots + A_k)^2 =  \sum A_i^2 + \sum_{i\ne j}
2A_i\cdot A_j \ge -1 - 2k + 2k.
$$
Hence all inequalities are equalities.  Thus the components are embedded and $A_i\cdot A_j$ is $0$ or $1$.  This means that the classes $A_i$ are all distinct.  But at most one class $A_i$ with $c_1(A_i)=0$ is represented for any $J_z$.  Thus $E$ has at most two components, and these intersect transversally 
because $A_0\cdot A_1=1$.

Since the energy $\om(A_i)$ is bounded above by $\om(E)$ and $J_z$ varies in a compact set, the usual compactness argument 
 implies that
there are only finitely many possibilities for the classes $A_i$.
Moreover, because the family $J_z$ is generic, each point $z$ at which a class $A_i$ with $c_1(A_i)=0$ is represented is isolated. The result follows.\end{proof}

\begin{rmk}\labell{rmk:e}\rm A similar argument shows that in a generic $3$-parameter family of $J_z$ no new nodal curves appear, though now of course degenerations appear for a $1$-parameter family of $z$.
\end{rmk}

The first (well known) corollary shows that in $4$ 
dimensions any form on a $k$-fold blow up that can be deformed  into
a blow up form must itself be a blow up form.

\begin{cor}  Consider a family $\Tom_t, t\in I,$ of symplectic forms on a $4$-dimensional manifold $\TM_k$ such that $\Tom_0$ is a blow up form.  Then $\Tom_1$ is a blow up form.
\end{cor}
\begin{proof} By assumption there are $k$ disjoint $\Tom_0$-symplectically embedded $-1$ spheres.  Denote their classes by $E_1,\dots, E_k$, and choose an $\Tom_0$-tame almost complex structure for which they are holomorphic.
Extend $J_0$ to a generic family $J_t, t\in I$, of 
almost complex structures such that $J_t$ is $\Tom_t$-tame.  By the lemma 
each class $E_i$ has a unique embedded $J_t$ holomorphic representative.  
Hence result.\end{proof}

%The next proposition implies 
We now show that many of
the new elements constructed in Proposition~\ref{prop:csympk2} 
do not lie in
$\pi_2(B\Symp(\TM_k,\Tom))$ for any blow up form. For simplicity, we
restrict to simply connected $M$ (so that every $M$-bundle over $S^2$
has a symplectic form) and will  
first formulate our argument in a special case.  We shall denote by  
$\Ee: = \Ee(M, \om)\subset H_2(M;\Z)$  the set
of classes
that are represented by symplectically embedded spheres of
self-intersection $-1$
(or symplectic $-1$ spheres, for short).
Recall from~\cite{Mrr} that the following conditions on a symplectic $4$-manifold are equivalent:
\begin{itemize}
\item $(M,\om)$ is not a blow up of $\C P^2$ or of a ruled surface (manifold with a fibering by symplectically embedded $2$-spheres);\SSS

\item
$E\cdot E'=0$ for any two distinct elements $E, E' \in \Ee.$  
\end{itemize}
   
If these conditions hold we shall say that $(M,\om)$ has Kodaira dimension $\ka(M)\ge0$.  Such manifolds have a unique minimal reduction; equivalently, if $C_1,\dots, C_\ell$ is a maximal set of disjoint symplectic $-1$ spheres then 
$\{[C_i]: 1\le i\le \ell\} = \Ee$.  Finally, $(M,\om)$ is said to be minimal if $\Ee(M) = \emptyset$.

\begin{prop}\labell{prop:kblow} Let $(M,\om)$ be a minimal symplectic $4$-manifold with $\ka(M)\ge 0$ and let 
$(\TM_k, \Tom)\to \TP \to S^2$ be a symplectic 
bundle with fiber the $k$-fold blow up of $(M,\om)$ with weights 
$\eps_1\ge\dots \ge \eps_k$.  Then $\TP$
 is the $k$-fold blow up of a symplectic 
 bundle $(M,\om) \to P \to S^2$
along sections $s_i, 1\le i\le k,$ whose homotopy
 class varies in a group spanned by $\pi_2(M;\Z)$ and the 
classes $E_j, j<i,$ of the exceptional divisors with $\eps_j>\eps_i.$
\end{prop}
\begin{proof}
Denote by $E_1,\dots, E_k\in H_2(\TM_k)$ the classes of the 
exceptional divisors; set $\eps_i: = \Tom(E_i)$. 
By the minimality of $(M,\om)$, the only classes in $\Ee(\TM_k)$
are the $E_i$. Let $J_z, z\in S^2$,
be a generic family of $\Tom_z$-tame 
almost complex structures on the fibers of $\TP\to S^2$, where
$\Tom_z$ is the symplectic form in the fiber $M_z: = \pi^{-1}(z)$ of $P\to S^2$.
Then we claim that the smallest class $E_k$ has an embedded $J_z$-holomorphic representative for every $z\in S^2$.  Otherwise, by 
an obvious generalization of Lemma~\ref{le:e}, we may write $E_k= A_0 +A_1$ where $A_0\in \Ee(\TM_k) = \{E_1,\dots,E_k\}$. Since $\Tom(A_1)< \Tom(E_k)=\eps_k
\le \eps_i$ for all $i$, this is impossible.

Let $\TSi$ be the submanifold of $\TP$ formed by the union of
these $E_1$-curves.  Then $\TSi$ is an $S^2$ bundle over $S^2$ with fibers that are $\Tom_z$-symplectically embedded.  
Choose a closed extension $\TOm$ of the $\Tom_z$ defined in some neighborhood $U$ of $\TSi$ and such that $\TSi$ is a symplectic submanifold of  $U$.\footnote
{
Note that we do not assume that $H^1(M;\R)=0$ or that
$\TP$ supports
 a closed extension of the $\Tom_z$, i.e. the bundle 
 $\TP\to S^2$ need not be Hamiltonian.  The present arguments  are geometric and do not use the coupling form.} 
  By ~\cite{Mrr} the restriction of $\TOm$ to $\TSi$ 
is  standard, and so, by the symplectic neighborhood theorem, $\TSi$ may be blown down.
Hence $\TP$ is formed by blowing up 
some bundle $\TQ\to S^2$ with fiber $\TM_{k-1}$ along some section $s_k$.  Repeating this argument, we find that $\TP$ is a $k$-fold blow up along sections $s_i$ whose homotopy class depends on $\pi_2(M)$ and 
the $E_j, j<i$.

It remains to prove the last statement. If $\eps_i >\eps_{i+1}$ for all $i$ there is nothing to prove.  Otherwise the result follows from Lemma~\ref{le:sk} below.
\end{proof}

We need to adapt Proposition~\ref{prop:size} to 
the case when the initial bundle $M\to P\to S^2$ is nontrivial.

\begin{lemma}\labell{le:sk}
Let $(X_1, \om_1)\to \TP_1\to S^2$ be 
 the $\eps_1$-blow up of a symplectic bundle 
$(X_0,\om_0)\to P\to S^2$ along a section $s_1$.  Denote 
the exceptional class in $X_1$ by $E_1$ and
 the bundle of exceptional divisors in $\TP_1$ by $\TSi_1$.
%
% Let $(X,\om)$ be a symplectic manifold that is itself a blow up of size $\al$, that is $\al$ is the integral of $\om$ over a line in the exceptional divisor $\Si_X$.   
Let $s_2:S^2\to \TP_1$ be a section such that $[s_2]\cdot [\TSi_1] \ne 0$, and construct    
$\pi: (\TP_2, \Tilde\Om_\eps)\to S^2$ as the $\eps_2$-blow
up of $\TP_1$ along $s_2$. Then
 $\eps_1 > \eps_2$.%\NI
%{\rm (ii)} If 
%$\Tilde\rho_t, 0\le t\le 1,$ is any family of symplectic forms on $\TP$ that start at some blow up form $\Tilde\rho_0: = \Tilde\Om_\eps$ then $\int_{E_\si}\Tilde\rho_1 < \al$,
%where $E_\si$ is the class of a line in
%the exceptional divisor in the blow up $\TX$.
\end{lemma}  
\begin{proof} The proof is essentially the same as that of Proposition~\ref{prop:size}.  We first show that
$$
\GW^{\TP_1,E_1}_{0,2}([s_2], [\TSi_1]) = - [s_2]\cdot [\TSi_1]\ne 0,
$$  
and then prove that 
$$
\GW^{\TP_1,E_1}_{0,2}([s_2], [\TSi_1]) = 
\GW^{\TP_2, E_1-E_2}_{0,2}([\TSi_2],[\TSi_1]),
$$  
where $\TSi_2$ is the bundle of exceptional divisors in $\TP_2$.
Further details are left to the reader.
\end{proof}

 It follows from Remark~\ref{rmk:e} that each blow down above is unique, i.e. it is independent of the choice of the family $J_z$.
 Moreover, it depends on the section only up to homotopy.
  Hence we find:
 
 \begin{cor}  With  $(\TM_k,\Tom_{\eps})$ as above, suppose 
 that $\eps_1>\dots>\eps_k>0$.
Then the rank
$ \pi_2\bigl(B\Symp(\TM_k,\Tom_{\eps})\bigr)$ exceeds that of
$\pi_2\bigl(B\Symp(M,\om)\bigr)$  by at most $rk+ k(k-1)/2$,
 where $r={\rm rk\,} \pi_2(M)$.
 \end{cor}

We cannot give a more precise answer here because
 there is a {\bf realization problem}: it is not obvious that the condition in Lemma~\ref{le:sk} is the only 
 obstruction to being able to blow up a section with weight $\eps$, even if we assume that $\ka(M)\ge 0$ and that $M$ itself has such a blow up.  Therefore some work is still needed before we can fully understand
 the relation between 
 $\pi_2\bigl(B\Symp(\TM_k,\Tom_{\eps})\bigr)$ and $\pi_2\bigl(B\Symp(M,\om)\bigr)$.\footnote
 {
 It is likely that by using appropriate local models one can show that
 given $(M,\om)$ and $k$ one can find $\eps_0$ so that
Lemma~\ref{le:sk} does 
  give the only obstructions when all $\eps_i\le \eps_0$.}
 Nevertheless, comparing with Proposition~\ref{prop:csympk2} it is clear that when $k\ge2$ there are many elements in 
 $\pi_2\bigl(B\Diff(\TM_k)\bigr)$ that are constructed using standard Hermitian structures but yet cannot be realized symplectically.
 Here is a sample result.
 
\begin{cor}\labell{cor:two}  Let $B_1, B_2$ be two disjoint embedded balls in a
symplectic $4$-manifold $(M,\om)$ with $\ka(M)\ge 0$. Denote by $X_i, i=1,2,$  the blow up of $M$ by $B_i$, and by 
$X_{12}$ the blow up by both balls.   Denote by $E_i$ the  class of the exceptional divisor
corresponding to $B_i$.  For $i=1,2$ let $\La_i\in\pi_2(B\Diff(X_{12}))$ be the element formed by blowing up $\TM_i\times S^2$ along the graph of $E_i$ with the standard Hermitian structure.
Then $\La_1 + \La_2$ is not homotopic to an element of $\pi_2\bigl(B\Symp(X_{12},\Tilde\om)\bigr)$  for any blow up form on $X_{12}$. 
\end{cor}
\begin{proof}  For $i = 1,2$ let  $X_{12}\to \TP_i\to S^2$ be the bundle corresponding to $\La_i$.  Then the bundle corresponding to
 $\La_1 + \La_2$ is the fiber\footnote
 {This is  the Gompf (or connect) sum of $\TP_1, \TP_2$ along a fiber.}
  sum $\TQ:\TP_1\#\TP_2\to S^2\# S^2 = S^2$; see~\cite{LMP}.
But $\TQ\to S^2$ is not a symplectic bundle since it does not satisfy the conclusion of Proposition~\ref{prop:kblow}.
\end{proof}

Finally we investigate the analog of Proposition~\ref{prop:kblow}
for blow ups of $\C P^2$. (We leave the case of 
 ruled surfaces to the reader.)

We need a preliminary  lemma.
Denote by $X_k$ the $k$-fold blow up of $\C P^2$, and fix a basis $L, E_1,\dots,E_k$ of $H_2(X_k;\Z)$ where $L$ is the class of a line and
$E_1,\dots E_k$ are represented by the obvious $-1$ spheres.
Let $\om^\mu$ be the standard (Fubini--Study) 
symplectic form on $X: = \C P^2$
normalized so that $\om^\mu(L)=\mu$.  
Define $\eps: = (\eps_1,\dots,\eps_k) $
where $\eps_1\ge \eps_2\ge \cdots$, and let $\om_\eps^\la$ be the symplectic form on $X_k$ obtained by cutting out $k$ disjoint balls from $(X,\om^\la)$ of sizes $\eps_1,\dots,\eps_k$ (assuming that this exists).  By \cite{Mcdef} all such forms are isotopic. We say that $\om_\eps^\la$ is 
{\bf reduced} if $\eps_1+\eps_2+\eps_3\le \mu$.

\begin{lemma}  Every symplectic form $\om$ on $X_k$ is diffeomorphic 
to some reduced form $\om_\de^\la$.
\end{lemma}
\begin{proof}  By \cite{LM2} $\om$ is diffeomorphic to some form $\om: = \om_\eps^\la$ with $\eps_1\ge \eps_2\cdots$ and so we just need 
to make sure we can arrange that $\eps_1+\eps_2+\eps_3\le\la.$
If $[\om]$ is rational then this is proved by B.-H.Li~\cite{BHL}.

His argument is as follows.  It suffices to suppose that  
$[\om]$  is integral.  Suppose that
$[\om]$ is not reduced and consider the form
$\om_{\eps'}^{\la'}$ obtained from $\om$ by doing the 
Cremona transformation in the first three $E_i$. This is a 
diffeomorphism $\phi_C$ that acts on homology by:
$$
\begin{array}{ccl}
L&\mapsto &L': = 2L-E_1-E_2-E_3,\\
E_i&\mapsto& E_i': = L-E_1-E_2-E_3+E_i,\;\; i=1,2,3\\
E_i&\mapsto& E_i':=E_i,\quad i>3.\end{array}
$$
Then $\la': = \om(L') < \la$, and the sequence $\eps'$
is obtained by  rearranging the numbers  $\om(E_i')$.
If $\om': = \om_{\eps'}^{\la'}$ is not reduced, then one repeats this process to obtain a new form with $\la''<\la'$.  
   Since these coefficients $\la', \la'',\dots$ 
    are positive integers, this process must stop after a finite number of steps. 

When $[\om]$ is not integral, we can argue as follows.  
Fix a generic complex structure $J$ on $X_k$, i.e. $(X_k,J)$ is the blow up of the complex manifold $\C P^2$ at $k$ generic points.
  It suffices to prove
 that 
 the numbers $\la', \la'', \dots$ lie in the finite set 
 $$
 S: = [0,\la]\cap
\{\om(A): A \mbox{ has a 
$J$-holomorphic representative}\}.
$$
Since there is a complex line in $\C P^2$ through $3$ generic points,
$L'$ has a $J$-holomorphic representative.  Hence 
$\la'\in S$.
But we can choose $\phi_C$ to be $J$-holomorphic.
Then $\la'':=\om'(L'') = (\phi_C^{-1})^*\om(L'') = 
\om((\phi_C^{-1})_*L'')$ lies in $S$ because the class $L''=
L'-E_1'-E_2'-E_3'$ and hence also $(\phi_C^{-1})_* L''$ 
have $J$-holomorphic representatives. Similar arguments apply to the subsequent numbers $\la''',\dots.$   \end{proof}

\begin{prop} Proposition~\ref{prop:kblow} and 
Corollary~\ref{cor:two} hold
when $M=\C P^2$.
% and the form $\om^\mu_\eps$ on its $k$-fold blow up $X_k$ is such that
%$\eps_1+\eps_2+\eps_3\le\mu$.
\end{prop}
\begin{proof}  
%By the previous lemma we may identify $(\TM_k,\Tom)$ with $(X_k,\om_{\eps})$ where $\eps_1+\eps_2+\eps_3<1$. 
By the previous lemma we may suppose that the form $\om_\eps^\mu$ on $\TM_k$ is reduced.  Hence, because 
 $\eps_i\ge \eps_{i+1}$ for all $i$, 
\begin{equation}\labell{eq:d}
\sum_{i\in I}\eps_i \le 3\mu\, d,
\end{equation}
for every
 $I\subset\{1,\dots,k\}$ with $3d$ elements.

 Let $\Ee$ be the set of classes in $H_2(X_{k})$ that can be represented by symplectically embedded $-1$ spheres.
 Then $\Ee$ contains the $E_i$, plus 
 many other classes such as $L-E_1-E_2$. However,
because every class in $\Ee$ has an embedded representative 
for generic $J$, positivity of intersections implies that every 
 class in $\Ee$ except the $E_i$ has the form
$dL - \sum m_iE_i$, where $m_i\ge 0$ and hence  $d> 0$.
To make the previous argument work we just need to show that 
in a generic $2$-dimensional family $J_z$
the class $E_k$ with minimal energy  never decomposes
into a nodal $J_z$-curve of type $(E, E_k-E)$, where $E\in \Ee$.

Suppose $E_k$ does decompose for $J_z$.  
By minimality,  $E\ne E_j$ and so we can write
 $E=dL - \sum m_iE_i$.
%It also satisfies the adjunction inequality
%$ A\cdot A - c_1(A) \ge -2$, i.e. $\sum m_i^2\le 1+d^2$.
Now observe that $E_k\cdot E = 0$ because $E\cdot(E_k-E) = 1$.
Hence $m_k=0$.  Let $F: = E_k-E$ be the $-2$-curve.
Choose $i$ so that $m_i>0$.
Then $E_i\cdot E=-m_i<0$ which means that $E_i$ also decomposes for $J_z$.  Since by genericity there is at most one $-2$ class with a $J_z$-holomorphic representative, $E_i$ must decompose as 
$(E',F)$ where $E'\in \Ee$ and $F$ is as above.  As before $E_i\cdot E' = 0$ and so $E_i\cdot F = m_i = 1$.
  Repeating this argument, we see that   $F = \sum m_iE_i-dL$, where $m_i$ is $0$ or $1$. Since $c_1(F) = 0$,
there must be precisely $3d$ nonzero $m_i$.  But then $\Tom(F)\le0$
by (\ref{eq:d}), which is impossible.   
\end{proof}

\begin{rmk}\rm For example, if $X_2$ is the $2$-fold blow up of $\CP^2$,  the above results show that $\pi_1(\Symp\, X_2)$ has rank at most $2$ when the blowups are of equal size and rank
at most $3$ otherwise.  In this case, the work of Lalonde and Pinsonnault~\cite{LP,Pin} shows that these are the precise ranks.\end{rmk}

The above results give examples of elements in $\pi_2(B\Diff)$
that are not themselves in $\pi_2(B\Symp)$ but are sums of such elements.  We now return to the question considered in 
Proposition~\ref{prop:diff}; are there elements of
$\pi_2(B\Diff)$ that are not in the subgroup generated by the 
different $\pi_2(B\Symp)$?
The next result illustrates what our current methods
can say about this.

\begin{prop}\labell{prop:diff2} 
Let $M$ be a simply connected minimal symplectic $4$-manifold 
with $\ka(M)\ge 0$ and a unique Seiberg--Witten basic class $K$. 
Let
 $\TM_k$ be its $k$-fold blow up for some $k>0$.   Then $\pi_2(B\Diff\TM_k)$
  is not generated by the images of  
  $\pi_2\bigl(B\Symp(\TM_k, \Tom)\bigr)$
 as $\Tom$ ranges over all symplectic forms on $\TM_k$.
 \end{prop}
\begin{proof} 
 T.-J. Li proved in~\cite[Cor.~3]{TJL} that 
when $\ka(M)\ge 0$ every symplectic form on $\TM_k$ is a blow up form with exceptional divisors in the same classes $E_1,\dots,E_k$.  In particular its
minimal reduction $M$ is unique up to diffeomorphism.
By work of Taubes~\cite{Tau}, the assumption that there is a unique
basic class implies that every symplectic structure on $M$ has the same homotopy class of almost complex structures.  Hence every symplectic
 bundle $(M,\om)\to P\to S^2$ has a well defined complex vertical tangent bundle whose first Chern class, denoted $c_1^V$, is independent of $\om$.\MS
 
% \NI
% {\bf Step 1:} {\it Choice of reference section.}\,\, 
% The next step is to choose a suitable reference section $s_0\in H^2(P;\R)$
% of a given symplectic bundle 
% $(M,\om)\to P\to S^2$.  If the homomorphisms $\pi_2(M)\to \R$ induced by $[\om]$ and $K=c_1(M)$ are linearly independent, (i.e. if $(M,\om)$ is {\it not} spherically monotone), then it is 
% easy to see that for every
%such bundle there is $s_0$ such that $I_a(s_0) = I_c(s_0) = 0$.  (For a proof, see ~\cite{LMP};  $I_a, I_c$ are defined in equations (\ref{eq:Ia}) and (\ref{eq:Ic}).) Note that when $I_c:\pi_2(M)\to \Z$ is nonzero, it suffices to consider $[\om]$ satisfying this condition, because if a class in $\pi_2(B\Diff\TM_k)$  is in the image of  $\pi_2\bigl(B\Symp(\TM_k, \Tom)\bigr)$ 
%it will remain in this image when $[\om]$ is slightly perturbed.
%Hence we may assume that the reference section $s_0$ satisfies 
%$$
%\mbox{ (i) } I_a(s_0) = I_c(s_0) = 0,\;\;\mbox{ or (ii) }
%I_a(s_0)=0,\; c_1(M)=0 \mbox{ on }\pi_2(M).
%$$
%In either case,  $s_0$ depends on the class $a: = [\om]$.\MS

 \NI
 {\bf Step 1:} {\it Volume calculations.}\,\,
By  Proposition~\ref{prop:kblow}, 
every symplectic bundle $(\TM_k,\Tom) \to (\TP, \TOm)\stackrel{\pi}\to S^2$ is the $k$-fold blow up of a 
possibly nontrivial bundle $(M,\om)\to P\stackrel{\pi}\to S^2$ along certain sections $s_i$
with weights $\eps_i$, where $\eps_1\ge \eps_2\ge\dots $.   
Thus there are  bundles 
$\TM_{i-1}\to \TP_i\to S^2$ for $ i=1,\dots, k$ with $\TP_k=\TP$.
Choose a reference section $s_0$ for  $P\to S^2$.  
Since we may assume this is disjoint from the $s_i$ there is a unique corresponding reference section, also called $s_0$, in each
bundle $\TP_i\to S^2$, and in particular for $\TP: = \TP_k$.

To calculate the volume of $P$,
choose a basis $B_\al, 1\le \al\le r,$ for $H_2(M;\Q)$
and dual basis $b_\al, 1\le \al\le r,$ for $H^2(M)$.  For each $a\in \Cc_M$, write 
$a= \sum \la_\al b_\al$.  Let $\Ha\in H^2(P)$ be its 
canonical extension and set 
$\Ha_\mu: = \Ha+\mu\pi^*(\be)$ where $\be\in H^2(S^2)$ has total area $1$.  Then $\Vol(P,\Ha_\mu)= \mu\,\Vol(M,a) = :\mu V_a$, where $V_a:=\Vol(M,a)$ is a polynomial function of the $\la_\al$.
%Every other section class $s$ of $P\to S^2$ may be written as $s_0+B$,
%where $B = \sum n_\al B_\al\in H_2(M)$.

%We now work out  the volume of $\TP$.  For any bundle
%$M\to P\to S^2$ the canonical extension $\Ha$ is defined so  that
%$\Vol(P,\Ha)=0$.  For each $\mu>0$ we define $\Ha_\mu: = \Ha+\pi^*(\be)$ where $\be\in H^2(S^2)$ has total area $\mu$.  It is easy to check that
%$\Vol(P,\Ha_\mu)= \mu\,\Vol(M,a) = :\mu V_a$, where $V_a:=\Vol(M,a)$ is a polynomial function of the $\la_\al$.  

To calculate the volume of $\TP$,
let $\Ha_{\mu,i,\eps}$ be the class on
$\TP_i$ obtained from $\Ha_\mu$ by blowing up along the sections
$s_j$ with weights $\eps_j$ for all $j\le i$. (Note that 
$\Ha_{\mu,i,\eps}$ is 
{\it not} a coupling class.)
Similarly, denote by $a_{i,\eps}$ its restriction to the fiber, i.e. the class on $\TM_i$ obtained by doing the first $i$ blow ups.
Then, by (\ref{eq:vol}),
%%%
%\begin{equation}\labell{eq:Tvol}
$$
\vol(\TP,\Ha_{\mu,k,\eps}) = \mu V_a -
\sum_i 
v_{\eps_i}\bigl(\int_{s_i} \Ha_{\mu,i-1,\eps}
-\frac{\ell_i}{n+1} \eps_i\bigr), \quad\mbox{where } \ell_i: = c_1^V(s_i).
$$
If we write $s_i = s_0 + B_i$ where $B_i = \sum n_{i\al} B_{\al}+\sum_{j<i} m_{ij}E_j$, then
%%%
%\begin{equation}\labell{eq:Tamu}
$$
\int_{s_i} \Ha_{\mu,i-1,\eps} = \mu + I_a(s_0) +\sum_\al n_{i\al}\la_\al
+\sum_{j< i}m_{ij}\eps_j,
$$
where $I_a$ is as in equation (\ref{eq:Ia}).  Thus $I_a(s_0)$
is a homogenous rational function of the $\la_\al$ of degree $1$.
Thus, 
\begin{equation}\labell{eq:Tvol}
\vol(\TP,\Ha_{\mu,k,\eps}) = \mu \bigl(V_a -
\sum_i v_{\eps_i}\bigr) - \sum_i%\Bigl( 
\Bigl(I_a(s_0) +\sum_\al 
n_{i\al}\la_\al +\sum_{j< i}m_{ij}\eps_j
-\frac{\ell_i}{n+1} \eps_i\Bigr)v_{\eps_i}.
\end{equation}
 Further,
since  this is a symplectic blow up, 
\begin{equation}\labell{eq:elli}
 \ell_i: = c_1^{V}(s_i) = I_c(s_0)+ \sum_\al n_{i\al} c_1(B_\al) + 
 \sum_{j< i}m_{ij}.
\end{equation}
Note that equation (\ref{eq:Tvol}) holds
 for all classes $a\in \Cc$, not just those represented by a symplectic form.
\MS

\NI
{\bf Step 2:} {\it The bundle $\TM_k\to\TQ\to S^2$}.
This bundle is obtained from the trivial bundle $(Q: = M\times S^2, \Om)$ by $k-1$ trivial blowups with weight $\eps_i$ 
and trivial Hermitian structure and with the last blow up
along $s_k: = gr_\si$ where $\si:S^1\to \TM_{k-1}$ represents $E_1$
and has a nonstandard Hermitian structure (i.e. $m_{k1}=1$ and $\ell_k\ne 1$.) 
Choose $s_0$ to be the trivial section.
Then, $I_a(s_0) = 0 = I_c(s_0)$ and, by (\ref{eq:Tvol}),
\begin{equation}\labell{eq:q}
\vol(\TQ,{\Ha_{\mu,k,\eps}}) = 
\bigl(V_a -\sum_{i} v_{\eps_i}\bigr)\mu - v_{\eps_k}\bigl(\eps_1 -
\frac{\ell_k}{n+1} \eps_k\bigr),
\end{equation}
for every $a\in \Cc$. 
\MS

\NI
{\bf Step 3:} {\it Sums of bundles.}
Since the sum in $\pi_2(B\Diff)$ corresponds to the fiber sum of bundles, it suffices to show that
$\TQ$ is not a finite fiber sum of symplectic bundles.  Suppose to the contrary that $\TQ = \TQ^{\ga_1}\#\dots\#\TQ^{\ga_m}  $
where for each $\ga: = \ga_i$ the bundle
 $\TQ^{\ga}\to S^2$ is the $k$-fold blow up of an 
$\om^\ga$-symplectic
$M$-bundle 
$Q^\ga\to S^2$. By the uniqueness of the blowdown, 
the sum $Q^{\ga_1}\#\dots\#Q^{\ga_m}$ is trivial.  The section $s_0$ 
splits into the sum $s_0^{\ga_1}\#\dots\#s_0^{\ga_m}$ of sections  of the $Q^\ga$. Note that this splitting is not unique.  However,
for all $a\in \Cc_M$,
$$
\sum_\ga I_a(s_0^\ga) = I_a(s_0) = 0,\quad 
\sum_\ga I_c(s_0^\ga) = I_c(s_0) = 0.
$$
As in Step 1, there are unique corresponding sections $s_0^\ga$ of the intermediate bundles $\TQ_i^\ga$.  Let $s_i^\ga$, $1\le i\le k$, be the section along which one does the $i$th blow up of $Q^\ga$.
Write $s_i^\ga = s_0^\ga + \sum n_{i\al}^\ga B_\al + 
\sum m_{ij}^\ga E_j$.\MS 
 
 \NI
{\bf Step 4:} {\it Completion of the proof.}
For each $\ga$ there is  a well defined class
$\Ha_{\mu^\ga}\in H^2(Q^\ga)$ that blows up to
$\Ha_{\mu^\ga,k,\eps}\in H^2(\TQ^\ga)$. Then, if $\mu: = \sum_\ga\mu^\ga$,
\begin{equation}\labell{eq:volkT}
\Vol(\TQ,{\Ha_{\mu,k,\eps}})= \sum_\ga \Vol(\TQ^\ga,{\Ha_{\mu^\ga,k,\eps_\ga}}).
\end{equation}
By (\ref{eq:Tvol}) the part of the coefficient  of $v_{\eps_i}$  in 
$
\sum_\ga \Vol(\TQ^\ga,{\Ha_{\mu^\ga,k,\eps}})$
that depends on the $\la_\al$ is 
$$
\sum_\ga I_a(s_0^\ga) + \sum_{\ga,\al} n_{i\al}^\ga \la_\al = 
\sum_{\ga,\al} n_{i\al}^\ga \la_\al.
$$
Comparing with (\ref{eq:volkT}) and (\ref{eq:q}), we find that 
$\sum_\ga n_{i\al}^\ga=0$ for all $i,\al$. Similarly, looking at the coefficient of $\nu_{\eps_i}\eps_j$ for $i\ne j$ we find  that
$\sum_\ga m_{ij}^\ga = 1$ if $(i,j)=(k,1)$ and $=0$ otherwise.
Because $\sum_\ga I_c(s_0^\ga) = 0$, it 
 now follows from equation (\ref{eq:elli}) that 
$\sum_\ga\ell_i^\ga = 1$ when $i=k$ and $=0$ otherwise.
Therefore by (\ref{eq:volkT}) the coefficient of $v_{\eps_k}\eps_k$
in $\Vol(\TQ,{\Ha_{\mu,k,\eps}})$ is $1/(n+1)$.
Since $\ell_k\ne 1$ this  contradicts equation (\ref{eq:q}).
\end{proof}

\begin{cor}\labell{cor:T4} Suppose that $\pi_1(M)\ne 0$ but that otherwise $M$ satisfies all the conditions in 
Proposition~\ref{prop:diff2}.  
Suppose further that $\cup a: H^1(M;\R)\to H^3(M;\R)$ is an isomorphism for all $a\in \Cc$. Then the conclusion of 
Proposition~\ref{prop:diff2} holds for the $k$-fold blow up $\TM_k$.
\end{cor}
 \begin{proof}   The condition $\pi_1(M)=0$ is used once in Step 4
  when we apply Lemma~\ref{le:coup} to assert 
  that every class $a\in \Cc_M$ extends to the total space of
 the symplectic $M$-bundles $Q^\ga\to S^2$, a necessary 
 step before we calculate volumes. 
  We now show that each of these bundles admits a section and 
  that therefore these extensions do exist.  The rest of the proof 
  then goes through as before.
  
  The fact that $Q^\ga\to S^2$ has a section is immediate from its construction as a blow down.   It remains to prove that the classes
  $a\in \Cc_M$ extend.  
%  
%The blow down argument of Proposition~\ref{prop:kblow} does not use 
%the closed extension $\TOm$ of the fiberwise symplectic forms $\om_z,
%z\in S^2$; all we need is a family $J_z, z\in S^2$, of $\om_z$-tame almost complex structures on the fibers of $\TQ^\ga\to S^2$.
%It follows from the existence of these blow downs that each of the intermediate bundles admits a section. In particular the initial bundle $(M,\om)\to Q^\ga\to S^2$ has a section.  We claim that 
%every class $a\in \Cc_M$ then extends to $Q^\ga$. 
%Since the base $S^2$ is simply connected, the normalization part of the proof of Lemma~\ref{le:coup} works as before.  
%This is enough for the proof of
%Proposition~\ref{prop:diff2} to go through.
The obstruction to extending the symplectic class $[\om]$ to the total space of a symplectic bundle $(M,\om)\to Q\to S^2$ is $\Flux^{[\om]}([\phi])$,
where the loop $\phi_t, t\in S^1,$ in $\Symp(M,\om)$  
is the clutching function of the bundle $Q=Q_\phi\to S^2$. 
Therefore the extension exists iff $\Flux^{[\om]}([\phi])=0$.
By our assumption on $M$, this is equivalent to requiring that the volume flux
$$
[\om]\wedge  \Flux^{[\om]}([\phi]) =: 
\Flux^{[\om]^2/2}([\phi])\in H^3(M;\R)
$$
 vanishes. 
The latter class is well known to be Poincar\'e dual to the image $\ev([\phi])$ of $[\phi]$ under the homology evaluation map
$\ev_*: \pi_1(\Diff(M))\to H_1(M)$. But 
 $Q_\phi\to S^2$ has a section  iff the loop
$t\mapsto \ev(\phi_t): = \phi_t(x_0)$ is contractible in $M$. 
Therefore, if there is a section 
$\ev_*([\phi]) = 0$.
% But, for any  $\si \in \pi_1\bigl(\Symp(M,\om)\bigr)$, 
%$ \Flux^{[\om]}(\si)\in H^1(M;\R)$ is zero  precisely if $[\om]\wedge  \Flux^{[\om]}(\si)\in H^3(M;\R)$ is zero, which happens precisely if the 
%Poincar\'e dual of this class vanishes in $H_1(M)$.  But this 
%Poincar\'e dual  is well known to be the image of
%$\ev(\phi_t)$.  
This shows that the symplectic class $[\om^\ga]$ extends to each $Q^\ga$.  Then every class $a\in \Cc_M$ has such an extension by ~\cite{LMP}. \end{proof}
 
 \begin{example}\rm  Taubes showed that 
every symplectic form on the $4$-torus $\T^4$ gives rise to an almost complex structure with $c_1=0$. Using this, or the fact that $\T^4$ 
has a unique Seiberg--Witten basic class (namely $K=0$),
 we can conclude that, when $M=\T^4$ and $k>0$, $\pi_2(B\Diff\TM_k)$
 is not generated by the images of  
$\pi_2\bigl(B\Symp(\TM_k, \Tom)\bigr)$
as $\Tom$ ranges over all symplectic forms on $\TM_k$.
\end{example}

 \begin{rmk}\rm   We have concentrated in this section on results on $\pi_2(B\Diff)$ since, by Lemma~\ref{le:e}, 
  Proposition~\ref{prop:kblow}
is valid for quite general choices of $[\om]$ and $\eps$
when the base is $2$-dimensional.  However, 
  if $[\om]$ is integral, there are special values for $\eps$
for which the blow up classes $E_i$ 
  never degenerate.  For example,
  one can either take the $\eps_i$ all equal to $1/N$,
 or, given positive integers $N_i$ 
 one could take $\eps_1=1/N_1,\eps_2=
\eps_1/N_2,
 \dots, \eps_k=  \eps_{k-1}/N_k$.   For such $\eps$ every symplectic bundle $(\TM_k,\Tom_k)\to \TP\to B$ is a $k$-fold blow up.  
This can be proved either by adapting  the proof of   
Proposition~\ref{prop:kblow} or by noting that the corresponding 
 spaces of (unparametrized) exceptional divisors 
$\Ee_i/U$ are all contractible.
 (Here $\Ee_i$ is the space of symplectically embedded $E_i$-spheres 
 in $\TM_{i-1}$; cf. the discussion at the beginning of \S2.)
In these cases our methods would
 give results about the higher homotopy and homology 
of $B\Symp\,\TM_k$ and $B\Diff\TM_k$.
\end{rmk}

\end{document}